\input amstex
\documentstyle{amsppt}

\def\gen{{\frak{g}}}

\def\a{{\alpha}}
\def\o{{\omega}}
\def\O{{\Omega}}
\def\G{{\Gamma}}
\def\l{{\lambda}}
\def\g{{\gamma}}
\def\G{{\Gamma}}
\def\b{{\beta}}
\def\eps{{\varepsilon}}

\def\1b{{\bold 1}}

\def\bb{{\bold b}}

\def\cb{{\bold c}}

\def\eb{{\bold e}}
\def\fb{{\bold f}}
\def\hb{{\bold h}}

\def\kb{{\bold k}}

\def\pb{{\bold p}}

\def\xb{{\bold x}}

\def\Ab{{\bold A}}
\def\Bb{{\bold B}}

\def\Ib{{\bold I}}

\def\Kb{{\bold K}}
\def\Nb{{\bold N}}

\def\Rb{{\bold R}}

\def\Ub{{\bold U}}

\def\Wb{{\bold W}}
\def\Xb{{\bold X}}

\def\Prod{{\ts\prod}}
\def\Wedge{{\ts\bigwedge}}

\def\tr{\roman{tr}}

\def\nor{{\text{nor}}}
\def\max{{{\text{max}}}}

\def\Hom{\text{Hom}\,}
\def\Id{\text{Id}\,}
\def\Proj{\text{Proj}\,}

\def\Det{\text{Det}\,}

\def\simto{\,{\buildrel\sim\over\to}\,}

\def\End{\text{End}\,}
\def\tr{\text{tr}\,}

\def\Ker{\text{Ker}\,}

\def\Im{\text{Im}\,}

\def\GL{{\text{GL}}}

\def\Oplus{\ts\bigoplus}
\def\Otimes{\ts\bigotimes}
\def\Sum{\ts\sum}

\def\AA{{\Bbb A}}

\def\CC{{\Bbb C}}
\def\DD{{\Bbb D}}

\def\NN{{\Bbb N}}

\def\QQ{{\Bbb Q}}

\def\ZZ{{\Bbb Z}}

\def\Bc{{\Cal B}}

\def\Ec{{\Cal E}}
\def\Fc{{\Cal F}}

\def\Oc{{\Cal O}}

\def\Qc{{\Cal Q}}

\def\Uc{{\Cal U}}
\def\Vc{{\Cal V}}
\def\Wc{{\Cal W}}

\def\and{{\quad\text{and}\quad}}

\def\ds{\displaystyle}                
\def\ts{\textstyle}

\def\qed{\hfill $\sqcap \hskip-6.5pt \sqcup$}        
\overfullrule=0pt                                    

\def\lra{{{\longrightarrow}}}

\def\iu{{\underline i}}
\def\ju{{\underline j}}

\def\u1{{\underline 1}}

\def\uj{{\underline j}}

\def\oh{{\overline h}}

\def\la{{\langle}}
\def\ra{{\rangle}}
\def\lla{{\longleftarrow}}

\newdimen\Squaresize\Squaresize=14pt
\newdimen\Thickness\Thickness=0.5pt
\def\Square#1{\hbox{\vrule width\Thickness
	      \alphaox to \Squaresize{\hrule height \Thickness\vss
	      \hbox to \Squaresize{\hss#1\hss}
	      \vss\hrule height\Thickness}
	      \unskip\vrule width \Thickness}
	      \kern-\Thickness}
\def\Vsquare#1{\alphaox{\Square{$#1$}}\kern-\Thickness}



\NoBlackBoxes

\topmatter

\title CANONICAL BASES AND QUIVER VARIETIES\endtitle
\author Michela Varagnolo and Eric Vasserot \endauthor

\address D\'epartement de math\'ematique, Universit\'e de Cergy-Pontoise,
2, av. A. Chauvin, BP 222, 95302 Cergy-Pontoise cedex, France\endaddress

\email michela.varagnolo\@math.u-cergy.fr\endemail

\address D\'epartement de math\'ematique, Universit\'e de Cergy-Pontoise,
2, av. A. Chauvin, BP 222, 95302 Cergy-Pontoise cedex, France\endaddress

\email eric.vasserot\@math.u-cergy.fr\endemail


\thanks
Both authors are partially supported by EU grant
\# ERB FMRX-CT97-0100.\endthanks


\abstract
We prove the existence of canonical bases in the $K$-theory of quiver varieties.
This existence was conjectured by Lusztig.
\endabstract
\endtopmatter
\document

\head Contents\endhead

1. Introduction

2. The algebra $\Ub$

3. The braid group

4. Reminder on quiver varieties

5. The involution on the convolution algebra

6. The metric and the involution on standard modules

7. Construction of the signed basis

8. Example

\head 1. Introduction\endhead
Lusztig proposed in \cite{16} to construct 
a signed basis of the equivariant $K$-theory of a quiver variety. 
As in \cite{14}, this signed basis should be characterized
by an involution and a metric.
He suggested a formula for the involution and the metric
and he conjectured the existence of the signed basis.
This signed basis should also satisfy some positivity property,
related, hopefully, to the positivity of the structural constants
of the product and the coproduct of the modified quantum algebra
in the canonical basis, for all simply laced types. 
The main purpose of this paper is to give a precise definition
of this signed basis and to prove its existence. 
It was conjectured in \cite{24} that
the $K$-theory of the quiver variety, 
with the action of the quantized enveloping algebra
of affine type defined in \cite{20}
(see also \cite{23} for the type $A$ case),
is isomorphic to the 'maximal integrable module'
introduced by Kashiwara in \cite{8}. 
This module has a canonical basis, see {\it loc. cit}. 
The conjectures in \cite{9, \S 13} suggest  
that Kashiwara's canonical basis and the geometric one are related,
see Remark 7.2.2.

We thank the referee for useful suggestions.

\vskip3mm

\head 2. The algebra $\Ub$\endhead
\subhead 2.1\endsubhead
Let $\gen$ be a simple, simply laced, complex Lie algebra.
Let $(a_{ij})_{i,j\in I}$ be the Cartan matrix.
The quantum loop algebra associated to ${\frak g}$ is the
$\QQ(q)$-algebra $\Ub'$ generated by 
$\xb^\pm_{ir},\,\kb^\pm_{is},\,\kb_i^{\pm 1}=\kb^\pm_{i0}$ 
$(i\in I,\,r\in\ZZ,\,s\in\pm\NN^\times)$
modulo the following defining relations 
$$\kb_i\kb_i^{-1}=1=\kb^{-1}_i\kb_i,
\quad [\kb^\pm_{i,\pm r},\kb^\eps_{j,\eps s}]=0,$$
$$\kb_i\xb^\pm_{jr}\kb_i^{-1}=q^{\pm a_{ij}}\xb^\pm_{jr},$$
$$(w-q^{\pm a_{ji}}z)\,\kb^\eps_j(w)\,\xb^\pm_i(z)=
(q^{\pm a_{ji}}w-z)\,\xb^\pm_i(z)\,\kb^\eps_j(w),$$
$$(z-q^{\pm a_{ij}}w)\xb^\pm_i(z)\xb^\pm_j(w)=
(q^{\pm a_{ij}}z-w)\xb^\pm_j(w)\xb^\pm_i(z),$$
$$[\xb^+_{ir},\xb^-_{js}]=\delta_{ij}{{\kb_{i,r+s}^+-\kb_{i,r+s}^-}
\over q-q^{-1}},$$
$$\sum_w\sum_{p=0}^m(-1)^{^p}
\left[\matrix m\cr p\endmatrix\right]
\xb^\pm_{ir_{w(1)}}\xb^\pm_{ir_{w(2)}}\cdots\xb^\pm_
{ir_{w(p)}}\xb^\pm_{js}\xb^\pm_{ir_{w(p+1)}}\cdots\xb^\pm_{ir_{w(m)}}=0,$$
where $i\neq j,$ $m=1-a_{ij},$ $r_1,...,r_m\in\ZZ,$ and $w\in S_m$.
We have set $[n]=q^{1-n}+q^{3-n}+...+q^{n-1}$ if $n\geq 0$,
$[n]!=[n][n-1]...[2]$, and
$$\left[\matrix m\cr p\endmatrix\right]=
{[m]!\over[p]![m-p]!}.$$
We have also set $\eps=+$ or $-$, and
$$\kb^\pm_i(z)=\sum_{r\geq 0}\kb^\pm_{i,\pm r}z^{\mp r},\quad
\xb_i^\pm(z)=\sum_{r\in\ZZ}\xb_{ir}^\pm\,z^{\mp r}.$$

\subhead 2.2\endsubhead
Put $\AA=\ZZ[q,q^{-1}]$.
Consider the $\AA$-subalgebra $\Ub\subset\Ub'$ 
generated by the quantum divided powers
$(\xb^\pm_{ir})^{(n)}=(\xb^\pm_{ir})^n/[n]!$, the Cartan elements
$\kb_i^{\pm 1}$, and the coefficients of the series
$$\sum_{s\geq 0}\pb_{i,\pm s}\, z^s=
\exp\Bigl(\sum_{s\geq 1}{\hb_{i,\pm s}\over [s]}z^s\Bigr),$$
where the elements $\hb_{is}$ are such that 
$$\kb^\pm_i(z)=\kb_i^{\pm 1}\exp
\Bigl(\pm(q-q^{-1})\sum_{s\geq 1}\hb_{i,\pm s}z^{\mp s}\Bigr).$$
Observe that $\Ub$ coincides with the $\AA$-subalgebra generated by the elements
$(\eb_i)^n/[n]!$, $(\fb_i)^n/[n]!$, and $\kb_i^{\pm 1}$, $i\in I\cup\{0\}$, 
where $\eb_i,\fb_i,\kb_i^{\pm 1}$
are the Kac-Moody generators, see \cite{3, Proposition 2.2 and 2.6}.

\subhead 2.3\endsubhead
Let $\Delta$ be the coproduct of $\Ub'$ defined in terms of the Kac-Moody 
generators as follows
$$\Delta(\eb_i)=\eb_i\otimes 1+\kb_i\otimes\eb_i,\quad 
\Delta(\fb_i)=\fb_i\otimes\kb_i^{-1}+1\otimes\fb_i,\quad
\Delta(\kb_i)=\kb_i\otimes\kb_i.$$
Let $\tau,\psi,S$ be the anti-automorphisms of $\Ub'$ such that
$$\tau(\eb_i)=\fb_i,\quad\tau(\fb_i)=\eb_i,
\quad\tau(\kb_i)=\kb^{-1}_i,\quad\tau(q)=q^{-1},$$ 
$$\psi(\eb_i)=q\kb_i\fb_i,\quad\psi(\fb_i)=q\kb_i^{-1}\eb_i,
\quad\psi(\kb_i)=\kb_i,\quad\psi(q)=q,$$ 
$$S(\eb_i)=-\eb_i\kb_i^{-1},\quad S(\fb_i)=-\kb_i\fb_i,
\quad S(\kb_i)=\kb^{-1}_i,\quad S(q)=q.$$ 
The map $S$ is the antipode.
Let $x\mapsto\bar x$ be the algebra automorphism of $\Ub'$ such that
$$\bar\eb_i=\eb_i,\quad\bar\fb_i=\fb_i,
\quad\bar\kb_i=\kb^{-1}_i,\quad\bar q=q^{-1}.$$ 

\subhead 2.4\endsubhead
Let $\dot\Ub'$ be the modified algebra of $\Ub'$,
and let $\dot\Ub$ be the corresponding $\AA$-form.
Let $\eta_\l\in\dot\Ub$ be the idempotent denoted by 
$1_\l$ in \cite{13, \S 23.1}.

\vskip3mm

\head 3. The braid group\endhead
\subhead 3.1\endsubhead
Let $P,Q,$ be the integral weight lattice, 
and the root lattice of $\gen$.
Let $\o_i,\a_i$, $i\in I$, 
be the fundamental weights and the simple roots.
Let $Q^+\subset Q$, $P^+\subset P$ be the subsemigroups generated
by the simple roots and the fundamental weights.
We set $\rho=\sum_{i\in I}\o_i$.
Let $a_i$, $i\in I$, be the positive integers such that the element
$\theta=\sum_{i\in I}c_i\a_i\in Q^+$ is the highest root.
The integer $\cb=1+\sum_ic_i$ is the Coxeter number of $\gen$.

Let $\delta$ be the smallest positive imaginary root
of the corresponding affine root system.
Recall that the affine root $\a_0$ is $\delta-\theta$.
We set $\hat P=P\oplus\ZZ\delta$.

Let $W$ be the Weyl group of $\gen$. 
Let $w_0\in W$ be the longest element.
The extended affine Weyl group is the semi-direct product
$\tilde W=W\ltimes P$.
For any element $w\in\tilde W$ let $l(w)$ be the length of $w$.
Let $s_i\in\tilde W$, $i\in I\cup\{0\}$, be the affine simple reflexions.
The affine Weyl group 
is the normal subgroup $\hat W\subset\tilde W$
generated by the elements
$s_i$, $i\in I\cup\{0\}$.
Let $\Gamma$ be the quotient group $\tilde W/\hat W.$ 
It is identified with a group of diagram automorphisms of the 
extended Dynkin diagram of $\gen$. 
In particular $\Gamma$ acts 
on $\Ub$, $\tilde W$ in the obvious way.

Let $B_W$, $B_{\tilde W}$ be the braid groups of $W$, $\tilde W$.
The group $B_{\tilde W}$ is generated by elements $T_w,$ $w\in\tilde W$, 
with the relation $T_wT_{w'}=T_{ww'}$ whenever $l(ww')=l(w)+l(w')$.
The group $B_W$ is the subgroup generated by the elements $T_w$, $w\in W$.
For simplicity we set $T_i=T_{s_i}$ for any $i\in I\cup\{0\}$, 
and $\theta_i=T_{\o_i}$ for any $i\in I$.
The group $B_{\tilde W}$ acts on $\Ub$ by algebra automorphisms. 
Let $T_i$ be the operator denoted by $T_{i,1}''$ in \cite{13, \S 37.1.3}. 
If $i\neq j$ we have
$$T_i(\eb_j)=\sum_{s=0}^{-a_{ij}}(-1)^{s}q^{-s}\eb_i^{(-a_{ij}-s)}\eb_j\eb_i^{(s)},
\quad T_i(\eb_i)=-\fb_i\kb_i,$$
$$T_i(\fb_j)=\sum_{s=0}^{-a_{ij}}(-1)^{s}q^{s}\fb_i^{(s)}\fb_j\fb_i^{(-a_{ij}-s)},
\quad T_i(\fb_i)=-\kb_i^{-1}\eb_i.$$ 
We have also $T_i(\kb_j)=\kb_j\kb_i^{-a_{ij}}$ for all $i,j$.

For a future use, we introduce the following notations : 

-- let $\sigma$ be the automorphism of $B_{\tilde W}$ 
such that $\sigma(T_w)=T^{-1}_{w^{-1}}$ for all $w\in\tilde W$;

-- for any $\a\in Q$, let $\Ub_\a\subset\Ub$ be the subset 
of the elements $x$ such that
$\kb_i\,x\,\kb_i^{-1}=q^{(\a,\a_i)}x$ for all $i$;

-- for any $i\in I$ let $\iu\in I$ be the unique element such that
$w_0(\a_i)=-\a_\iu$;

-- let $(,)\,:\,P\times P\to\QQ$ be the pairing such that 
$(\o_i,\a_j)=\delta_{ij}$; 

-- for any $\a\in Q$ we set $|\a|^2=\sum_i(\o_i,\a)^2$.

\subhead 3.2\endsubhead
Let $\g s_{i_1}s_{i_2}\cdots s_{i_k}$ be a reduced expression for 
the element $\o_i\in\tilde W$. 
Set 
$$\g_i=\sum_{\ell=1}^k\gamma s_{i_1}\cdots s_{i_{\ell-1}}(\a_{i_\ell})
\in\hat P.$$

\proclaim{Lemma}
We have $(\g_i,\a_i)=-\cb$.
\endproclaim

\demo{Proof}
Let $\Delta_\pm\subset\bigoplus_{i\in I}\ZZ\a_i$ be the sets of positive and
negative roots.
Let $\hat\Delta_\pm\subset\Delta+\ZZ\delta$ be the sets of positive and
negative affine roots.
We put $\Delta=\Delta_+\sqcup\Delta_-$,
$\hat\Delta=\hat\Delta_+\sqcup\hat\Delta_-$ and
$\hat\Delta(\o_i)=\hat\Delta_+\cap\o_i(\hat\Delta_-).$
Then,
$$\g_i=\sum_{\b\in\hat\Delta(\o_i)}\b.$$
Recall that
$$\hat\Delta_+=\Delta_+\cup\bigcup_{n\geq 1}(n\delta+\Delta),\quad
\hat\Delta_-=\Delta_-\cup\bigcup_{n\geq 1}(-n\delta+\Delta),$$
and that $\o_i(\a)=\a-(\o_i,\a)\delta$ for all affine root $\a$.
Thus,
$$\hat\Delta(\o_i)=\{\a-(n-a_i)\delta\,|\,\a\in\Delta_-,\,a_i>n\geq 0\},$$
where we set $a_i=-(\o_i,\a)$. Thus,
$$\g_i=\sum_{\a\in\Delta_-}a_i\bigl(\a+{1+a_i\over 2}\delta\bigr).$$
Let $\kappa$ be the Killing form.
We get
$$\matrix
(\g_i,\a_i)
&=-\sum_{\a\in\Delta_+}(\o_i,\a)\cdot(\a_i,\a)\hfill\cr
&=-\kappa(\o_i,\a_i)/2\hfill\cr
&=-\cb,\hfill
\endmatrix$$
see \cite{6, Exercice 6.2}.
\quad\qed
\enddemo

We fix the Drinfeld generators of $\Ub$ in such a way that
$$\xb^-_{ir}=o_i^r\theta_i^r(\fb_i), \quad
\xb^+_{ir}=o_i^r\theta_i^{-r}(\eb_i),\leqno(3.2.1)$$
where $o_i=\pm 1$ and $o_i+o_j=0$ if $a_{ij}<0$, see \cite{2, Definition 4.6}. 
Note that there is exactly two choices for the map $i\mapsto o_i$.
A case-by-case computation shows that the 
integer $o_io_\iu$ does not depend on $i$ : 
it is equal to $(-1)^\cb$.

\proclaim{Proposition}
\roster
\item There are unique $\AA$-algebra automorphisms $A,B\,:\,\Ub\to\Ub$ such that
$$A(\xb_{ir}^\pm)=-q^{\mp 1}\xb_{ir}^\pm,\quad 
B(\xb_{ir}^+)=-\xb_{ir}^+\kb_i,\quad 
B(\xb_{ir}^-)=-\kb_i^{-1}\xb_{ir}^-.$$

\item We have $\tau(\xb^\pm_{ir})=\xb^\mp_{i,-r},$
$\tau(\kb^\pm_{ir})=\kb^\mp_{i,-r}.$

\item We have $\psi(\xb_{ir}^\pm)=q^{-r\cb}T_{w_0}A(\xb^\pm_{\iu,-r}),$
$\overline{\xb_{ir}^\pm}=q^{r\cb}T_{w_0}B(\xb^\mp_{\iu r}).$
\endroster
\endproclaim

\demo{Proof}
Claim 2 is known, see \cite{2}.
Claim 1 is a consequence of the identities 3.
Let us prove 3.
Let $\overline T_i$, ${}^\psi T_i$ be the automorphisms 
of the algebra $\Ub$ such that
$\psi(T_i(x))={}^\psi T_i(\psi(x))$,
$\overline T_i(x)=\overline{T_i(\overline x)}$ for all $x\in\Ub.$
By \cite{13, \S 37} we have $\overline T_i=T_{i,-1}''$.
A case-by-case computation gives also ${}^\psi T_i=T_{i,-1}''$.
If $x\in\Ub_\a$, $\a\in Q$, we have
$$T_{i,-1}''(x)=(-q)^{-(\a,\a_i)}T_i^{-1}(x)$$
for all $i$, see \cite{13, \S 37}. Thus,
$$\matrix
{}^\psi\theta_i(x)&=\overline\theta_i(x)\hfill\cr
&=(-q)^{-(\a,\b_i)}\sigma(\theta_i)(x)\hfill\cr
&=(-q)^{(\a,\g_i)}\sigma(\theta_i)(x),\hfill\cr
\endmatrix$$
where 
$\b_i=\a_{i_k}+s_{i_k}(\a_{i_{k-1}})+\cdots+s_{i_k}\cdots s_{i_2}(\a_{i_1})$.
Note that $(\a,\g_i)=-(\a,\b_i)$ since $\g_i=-\o_i(\b_i)$.
The weight $\o_i$ beeing dominant we have
$T_{w_0}T_{\o_i}=T_{-\o_\iu}T_{w_0}$, i.e. 
$T_{w_0}\theta_i T_{w_0}^{-1}=\sigma(\theta_\iu)^{-1}$. 
Recall that 
$$T_{w_0}(\eb_i)=-\fb_\iu\kb_\iu,\quad
T_{w_0}(\fb_i)=-\kb_\iu^{-1}\eb_\iu,\quad
T_{w_0}(\kb_i)=\kb_\iu^{-1},\quad\forall i\neq 0.$$
Note that $\theta_i(\kb_i)=\kb_i$, see \cite{2}.
Using (3.2.1) we get
$$\matrix
\overline{(\xb_{ir}^+)}
&=o_i^r(-q)^{r\cb}\sigma(\theta_i)^{-r}(\eb_i)\hfill\cr
&=-o_i^r(-q)^{r\cb}T_{w_0}\theta_\iu^{r}(\kb_\iu^{-1}\fb_\iu)\hfill\cr
&=-q^{r\cb}T_{w_0}(\kb_\iu^{-1}\xb^-_{\iu r}).\hfill
\endmatrix$$
Similarly we have
$$\matrix
\psi(\xb_{ir}^+)
&=o_i^r(-q)^{-r\cb}\sigma(\theta_i)^{-r}(q\kb_i\fb_i)\hfill\cr
&=-o_i^r(-q)^{-r\cb}T_{w_0}\theta_\iu^r(q^{-1}\eb_\iu)\hfill\cr
&=-q^{-1-r\cb}T_{w_0}(\xb^+_{\iu,-r}).\hfill
\endmatrix$$
The case of $\xb_{ir}^-$ is identical.
\quad\qed
\enddemo

\head 4. Reminder on quiver varieties\endhead
\subhead 4.1\endsubhead
Let the couple $(J,H)$ denote the quiver such that 
$J$ is the set of vertices,
$H$ is the set of arrows.
If $h\in H$ let $h',h''\in J$ be
the incoming and the outcoming vertex of $h$. 
Let $\oh$ denote the arrow opposite to $h$.
We will consider the following cases : 

- $\Pi=(I,H)$ where $I$ is as in 2.1 and $H$
is such that there are $2\delta_{ij}-a_{ij}$ arrows
from $i$ to $j$ for all $i,j$.
Then, let $\O\subset H$ be any set such that
$H=\Omega\sqcup\bar\Omega$. 
Let $n_{ij}$ (resp. $\bar n_{ij}$)
be the number of arrows in $\O$ (resp. $\bar\O$)
from $i$ to $j$. Note that $n_{ij}=\bar n_{ji}$.

- Fix a set $I^1$ with a bijection $I\simto I^1$, $i\mapsto i^1$.
The quiver $\Pi^e=(I^e,H^e)$ is such that
$I^e=I\sqcup I^1,$
$H^e=H\sqcup\{i\to i^1, i^1\to i\,|\,i\in I\}.$

\subhead 4.2\endsubhead
Fix $V=\bigoplus_{i\in I}\CC^{a_i}$,
$W=\bigoplus_{i\in I}\CC^{\ell_i}$.

\vskip3mm

\noindent{\bf Convention.}
Fix $(m_i)\in\ZZ^I$.
Hereafter let $\mu,\l,\a$ denote elements
in $P$, $P^+$, $Q^+$ respectively such that
$\mu=\sum_im_i\o_i$, 
$\l=\sum_i\ell_i\o_i$, 
$\a=\sum_ia_i\a_i$. 
The dimension of the graded vector space $V$
is identified with the root $\a$ 
while the dimension of $W$ is 
identified with the weight $\l$.

\vskip3mm

\noindent
The space
$$M_{\l\a}=\bigoplus_{h\in H}M_{a_{h'}a_{h''}}(\CC)\oplus
\bigoplus_{i\in I}\bigl(M_{a_i\ell_i}(\CC)\oplus
M_{\ell_i a_i}(\CC)\bigr)$$
is identified with the set of representations of the quiver $\Pi^e$
on $V\oplus W$.
For any $(B,p,q)\in M_{\l\a}$ let $B_h$ be the component of the element
$B\in\Hom(V_{h''},V_{h'})$ and set
$$m_{\l\a}(B,p,q)=\sum_{h\in H}\eps(h)B_hB_\oh +pq\in\bigoplus_i\Hom(V_i,V_i),$$
where $\eps$ is a function $\eps\,:\, H\to\CC^\times$ such
that $\eps(h)+ \eps(\oh)=0.$ 
Put $G_\l=\prod_i\GL_{\ell_i}$, $G_\a=\prod_i\GL_{a_i}$.
The group $\CC^\times\times G_\l\times G_\a$ acts on $M_{\l\a}$ by 
$$(z,g_\l,g_\a)\cdot (B,p,q)=
(z g_\a Bg_\a^{-1},zg_\a pg_\l^{-1},z g_\l qg_\a^{-1}).$$
Following \cite{19} , we consider the varieties 
$$Q_{\l\a}^{(\mu)}=\Proj\bigl(\Oplus_{n\geq 0}A_n^{(\mu)}\bigr)\and 
N_{\l\a}=m^{-1}_{\l\a}(0)/\!\!/G_\a,$$
where $/\!\!/$ is the categorical quotient, 
$$A_n^{(\mu)}=\bigl\{f\in\CC[m^{-1}_{\l\a}(0)]\,\big|
\,f\bigl(g_\a\cdot(B,p,q)\bigr)=\chi_\mu(g_\a)^{-n}f(B,p,q)\bigr\},$$ 
and $\chi_\mu(g_\a)=\Prod_i\Det(g_{a_i})^{m_i}.$ 
The obvious projection $\pi_{\l\a}\,:\, Q^{(\mu)}_{\l\a}\to N_{\l\a}$ 
is a projective map. 
If $\mu,\mu'$ are such that $m_i, m'_i>0$ for all $i$,
or $m_i, m'_i<0$ for all $i$,
then the varieties $Q^{(\mu)}_{\l\a}$, $Q_{\l\a}^{(\mu')}$ 
are canonically isomorphic.
There is an open subset
$m_{\l\a}^{-1}(0)^{(\mu)}\subset m_{\l\a}^{-1}(0)$
whose points are called $\mu$-semistable,
such that there is a good quotient of
$m_{\l\a}^{-1}(0)^{(\mu)}$ by the group $G_\a$ and we have
$$m_{\l\a}^{-1}(0)^{(\mu)}/\!\!/G_\a=Q_{\l\a}^{(\mu)},\leqno(4.2.1)$$
see \cite{18, \S 1.7} for instance.
Moreover, if $\mu$ is a regular weight
then (4.2.1) is a geometric quotient
and the variety $Q_{\l\a}^{(\mu)}$ is smooth,
see \cite{21, Proposition 2.6}.
If $\mu$ is regular dominant, i.e. if $m_i>0$ for all $i$,
we set $Q_{\l\a}=Q_{\l\a}^{(\mu)}$.

\vskip3mm

\noindent{\bf Convention.}
Hereafter, we assume that $(\mu,\a)\neq 0$ for any root $\a$.

\subhead 4.3\endsubhead
Put $d_{\l\a}=\dim Q_{\l\a}$. It is known that $d_{\l\a}=(\a,2\l-\a)$.
If $\a\geq\b$ the extension by zero of representations of the quiver
gives a closed embedding $N_{\l\b}\hookrightarrow N_{\l\a}$.
For any $\a,\a',$ we consider the fiber product
$$Z_{\l\a\a'}=Q_{\l\a}\times_\pi Q_{\l\a'}.$$
If $\a'=\a+n\a_i$, $n>0$,
let $X_{\l\a\a'}\subset Z_{\l\a\a'}$
be the set of pairs $(x,x')$ 
which are the $G_\a$-orbits of $\Pi^e$-modules $y,y'$ in
$M_{\l\a}$, $M_{\l\a'}$ with $y$ a subrepresentation of $y'$.
If $\a'=\a-n\a_i$, put
$X_{\l\a\a'}=\phi\bigl(X_{\l\a'\a}\bigr)\subset Z_{\l\a\a'}$,
where $\phi$ is the automorphism of $Q_\l\times Q_\l$
taking an element $(x,y)$ to $(y,x)$.
The variety $X_{\l\a\a'}$ is smooth, see \cite{20, \S 5.3}. 
Consider the following varieties 
$$N_\l=\bigcup_\a N_{\l\a},\quad Q_\l=\bigsqcup_\a Q_{\l\a},\quad
Z_\l=\bigsqcup_{\a,\a'}Z_{\l\a\a'},\quad
X_\l=\bigsqcup_{\a,\a'}X_{\l\a\a'},\quad F_\l=\bigsqcup_\a F_{\l\a},$$
where $\a,\a'$ take all the possible values in $Q^+$
and $F_{\l\a}=\pi_{\l\a}^{-1}(0)$.

\subhead 4.4\endsubhead
For any complex algebraic linear group $G$, and
any quasi-projective $G$-variety $X$ let $\Kb^G(X)$
be the Grothendieck group of $G$-equivariant 
coherent sheaves on $X$.
We put $\Rb^G=\Kb^G(point)$.
Let $\Xb^G\subset\Rb^G$ be the set of the simple modules.
If the $G$-equivariant sheaf $\Ec$ is locally free,
let $\wedge^i\Ec$ is its $i$-th wedge power, 
and $\Wedge_\Ec$ be its maximal wedge power.
Note that $\Wedge_\Ec$ is still defined, in the obvious way,
whenever $\Ec$ is a $G$-equivariant complex on $X$.

\vskip3mm

\noindent{\bf Convention.}
Hereafter, let $f_*, f^*,\otimes,$ denote the derived functors 
$Rf_*, Lf^*,\otimes^L$ when they exist. 
Here $\otimes$ is the tensor product of coherent sheaves.
We use the same notation for a sheaf 
and its class in the Grothendieck group. 

\subhead 4.5\endsubhead
Set $\tilde G_\l=G_\l\times\CC^\times$.
Let $q$ denote also the character of the group 
$\CC^\times$ such that $z\mapsto z$. 
The canonical bundle of the variety $Q_{\l\a}$ is
$$\O_{Q_{\l\a}}=q^{-d_{\l\a}},\leqno(4.5.1)$$
see \cite{24, \S 6.4} for instance.
Let $V_i,W_i$ be the vectorial representations of the groups
$\GL_{a_i}$, $\GL_{\ell_i}$. 
Consider the following elements in $\Rb^{\tilde G_\l\times G_\a}$
$$F_i^+=q^{-1}W_i-q^{-2}V_i+q^{-1}\sum_{a_{ij}=-1}V_j,\quad
F_i^-=-V_i,\quad
F_i=F_i^++F_i^-.$$
The group $\tilde G_\l$ acts on the variety $Q^{(\mu)}_{\l\a}$.
If $E$ is a $\tilde G_\l\times G_\a$-module, 
let $E^{(\mu)}=m_{\l\a}^{-1}(0)^{(\mu)}\times_{G_\a}E$ 
be the induced $\tilde G_\l$-bundle on $Q^{(\mu)}_{\l\a}.$ 
There is a unique ring homomorphism
$$\Rb^{\tilde G_\l\times G_\a}\to\Kb^{\tilde G_\l}(Q^{(\mu)}_{\l\a})$$
such that $E\mapsto E^{(\mu)}$ for all $E$.
If $\mu$ is dominant
we set $\Vc_i=V_i^{(\mu)}$
and similarly for $\Wc_i$, $\Fc_i^\pm$, $\Fc_i$.
We set also
$\Vc=\bigoplus_i\Vc_i$, $\Wc=\bigoplus_i\Wc_i$.

\vskip3mm

\noindent{\bf Convention.}
The restriction to $Q_{\l\a}^{(\mu)}$
of a sheaf $\Ec$ on $Q^{(\mu)}_\l$ is denoted by $\Ec_\a$.
For simplicity we set $\Fc_{i;\a}=(\Fc_i)_\a$, etc.

\vskip3mm

\subhead 4.6\endsubhead
Consider the map
$$\dagger\,:\,M_{\l\a}\to M_{\l\a},\ 
(B,p,q)\mapsto(B,p,q)^\dagger=(-\eps {}^tB,-{}^tq,{}^tp),$$
where the upperscript $t$ stands for the transpose map.
Note that $\dagger$ does not commute 
to the action of the group $\tilde G_\l\times G_\a$. 
Let $\dagger$ be the group automorphism of $\tilde G_\l\times G_\a$ such that
$$\dagger\,:\,(z,g_\l,g_\a)\mapsto(z,g_\l,g_\a)^\dagger=
(z,{}^tg_\l^{-1},{}^tg_\a^{-1}).$$
Then $(g\cdot x)^\dagger=g^\dagger\cdot x^\dagger$ for all 
$g\in\tilde G_\l\times G_\a$, $x\in M_{\l\a}$.
The induced map $\dagger\,:\,Q^{(\mu)}_{\l\a}\to Q^{(-\mu)}_{\l\a}$ 
is an isomorphism of algebraic varieties. 
Let 
$$\dagger\,:\,
\Rb^{\tilde G_\l\times G_\a}\to\Rb^{\tilde G_\l\times G_\a},
\,E\mapsto E^\dagger$$
be the ring automorphism 
induced by the group automorphism $\dagger$.
For any element $E\in\Rb^{\tilde G_\l\times G_\a}$,
let
$(E^{(\mu)})^\dagger\in\Kb^{\tilde G_\l}(Q^{(\mu)}_{\l\a})$ 
be the pull-back of 
$E^{(-\mu)}\in\Kb^{\tilde G_\l}(Q^{(-\mu)}_{\l\a})$ 
by the automorphism $\dagger$.
We have
$$(E^{(\mu)})^\dagger=(E^\dagger)^{(\mu)}.\leqno(4.6.1)$$

For any $w\in W$ we set $w*\a=\l-w(\l)+w(\a).$
The element $w*\a$ depends on the weight $\l$.
However, since $\l$ is fixed in the whole paper the 
notation $w*\a$ should not make any confusion.
There is a $\tilde G_\l$-equivariant isomorphism of algebraic varieties
$S_w\,:\,Q_{\l\a}^{(\mu)}\to Q_{\l,w*\a}^{(w(\mu))}$ for each $w$,
such that 
$$S_{i}^2=1\and S_{ww'}=S_wS_{w'}\ \text{if}\ l(ww')=l(w)+l(w'),$$
see \cite{15}, \cite{18}, \cite{21}
(for simplicity we set $S_i=S_{s_i}$,
where $s_i$ is the simple reflexion with respect to the root $\a_i$). 
The precise definition of $S_w$ is given in the proof of Lemma 4.6.
Consider the composed map $\o=S_{w_0}\,\dagger.$ 
This choice is motivated by \cite{16} and \cite{21, Theorem 11.7}.
The map $\o$ is an isomorphism of algebraic varieties $Q_{\l\a}\simto Q_{\l,w_0*\a}$.

\proclaim{Lemma} 
\roster
\item We have 
$\o^*(\Fc_i)=-q^{\cb}\Fc^\dagger_\iu,$
$\o^*(\Wc_i)=\Wc_i^\dagger$
and
$$\sum_j[a_{ij}]\bigl(\o^*(\Vc_j)+q^{\cb}\Vc_\ju^\dagger\bigr)=
\Wc_i^\dagger+q^{\cb}\Wc_\iu^\dagger.$$

\item We have $\o^2=\Id$.

\item We have $\o(F_{\l\a})=F_{\l,w_0*\a}$.
\endroster
\endproclaim

\demo{Proof}
We use the construction of the operator $S_w$ given in \cite{18}, 
see also \cite{15}. Let us recall it briefly.
Set $\a'=s_i*\a$, $\mu'=s_i(\mu)$, $\mu'=\sum_im'_i\o_i$.
Let first assume that $m_i<0$.
Then $m'_i>0$.
Following \cite{15, \S 3.2}, let 
$$Z_i^\mu\subset 
m_{\l\a}^{-1}(0)^{(\mu)}\times m_{\l\a'}^{-1}(0)^{(\mu')}$$
be the set of pairs $(x,x')$, 
where $x=(B,p,q)$, $x'=(B',p',q')$ are such that

-- the sequence of $\tilde G_\l\times G_\a$-modules 
$$0\,\lra\, q^{-2}V'_i{\buildrel a(x')\over\lra}\, 
q^{-1}W_i\oplus q^{-1}\bigoplus_{a_{ij}=-1}V_j
\,{\buildrel b(x)\over\lra}\, V_i\,\lra\, 0$$
such that 
$a(x')=(q'_i,B'_h)$, $b(x)=p_i+\eps_{\bar h}B_h$ is exact, 

-- we have $a(x)b(x)-a(x')b(x')=0$, 

-- we have $B_h=B'_h$ if $h',h''\neq i$, and $p_j=p'_j$, $q_j=q'_j$ if $j\neq i$. 

\noindent
Note that $a(x')$ is injective and $b(x)$ is surjective,
see \cite{18, Lemma 38}.
Thus $Z_i^\mu$ is a closed subset of
$m_{\l\a}^{-1}(0)^{(\mu)}\times m_{\l\a'}^{-1}(0)^{(\mu')}$.
Consider the group 
$G_{\a\a'}=GL_{a_i}\times GL_{a'_i}\times\prod_{j\neq i}GL_{a_j}$.
The categorical quotient 
$$Q_i^\mu= Z_i^\mu /\!\!/G_{\a\a'},$$
is a smooth variety.
Moreover, the obvious projections are
isomorphisms of algebraic varieties
$$p_{1,\a}^{(\mu)}\,:\,Q_i^\mu\simto\, Q^{(\mu)}_{\l\a},\quad
p_{2,\a'}^{(\mu')}\,:\,Q_i^\mu\simto\,Q^{(\mu')}_{\l\a'},$$
see \cite{18, Proposition 40}.
The group $\tilde G_\l$ acts in the obvious way on $Q_i^\mu$, 
making the maps $p_{1,\a}^{(\mu)}$, $p_{2,\a'}^{(\mu')}$ equivariant.
By construction, for any $i\neq j$ we have 
$$\matrix
(p_{1\a}^{(\mu)})^*(F_i^{(\mu)}+q^{-2}V_i^{(\mu)})=
(p_{2\a'}^{(\mu')})^*(q^{-2}V_i^{(\mu')}),\hfill\cr\cr
(p_{1\a}^{(\mu)})^*(V_j^{(\mu)})=
(p_{2\a'}^{(\mu')})^*(V_j^{(\mu')}).\hfill
\endmatrix\leqno(4.6.2)$$
We set (recall that $m_i<0$)
$$S_i=p_{2\a'}^{(\mu')}({p_{1\a}^{(\mu)}})^{-1}
\,:\,Q^{(\mu)}_{\l\a}\to Q^{(\mu')}_{\l\a'}.$$
If $m_i>0$ we set
$$S_i=p_{1\a'}^{(\mu')}({p_{2\a}^{(\mu)}})^{-1}
\,:\,Q^{(\mu)}_{\l\a}\to Q^{(\mu')}_{\l\a'}.$$
Using (4.6.2) we get, if $m_i<0$ and $i\neq j$,
$$S_i^*(V_i^{(\mu')})=q^2F_i^{(\mu)}+V_i^{(\mu)},\quad
S_i^*(V_j^{(\mu')})=V_j^{(\mu)}.$$
Note that the map $S_i$ commutes to the action of the group $\tilde G_\l$. 
Thus, 
$$S_i^*(W_j^{(\mu)})=W_j^{(\mu')}$$
for all $j$.
Set $\eps_i=+1$ if $m_i>0$, $\eps_i=-1$ if $m_i<0$.
A case-by-case analysis gives
the following equalities in $\Kb^{\tilde G_\l}(Q_{\l}^{(\mu')})$
$$S_i^*(F_j^{(\mu)})=
\left\{\matrix
-q^{2\eps_i}F_j^{(\mu')}\hfill&\text{if}\quad i=j\hfill\cr
F_j^{(\mu')}\quad\hfill&\text{if}\quad a_{ij}=0\hfill\cr
F_j^{(\mu')}+q^{\eps_i}F_i^{(\mu')}\quad\hfill&\text{if}\quad a_{ij}=-1.
\hfill
\endmatrix\right.$$
The general formula is 
$$S_i^*(F_j^{(\mu)})=
F_j^{(\mu')}-q^{\eps_i}[a_{ij}] F_i^{(\mu')}.\leqno(4.6.3)$$
We now assume that the weight $\mu$ is dominant.
Thus, $F_j^{(\mu)}=\Fc_j$.
Fix an element $w$ in the Weyl group.
Let us prove that 
$$w(\a_i)=\a_j\,\Rightarrow
S_{w}^*(\Fc_j)=q^{a(w,i)}F_i^{(w^{-1}(\mu))},\leqno(4.6.4)$$
where 
$$a(w,i)={1\over 2}\sum_{\a\in\Delta_+\cap w^{-1}\Delta_-}(\a_i,\a)^2.$$
We may assume that $l(w)>0$ and that (4.6.4) holds for 
any $x$ with $l(x)<l(w)$. 
Fix $k\in I$ such that $w(\a_{k})\in -Q^+$. 
Let $\la s_i,s_{k}\ra$ be the subgroup generated by $s_i$, $s_{k}$. 
Let $x$ be the element of minimal length in the set $w\la s_i,s_{k}\ra$. 
Then, $x(\a_i), x(\a_{k})\in Q^+$. 
One of the following two cases holds, see \cite{12, Proof of Proposition 1.8}.

\vskip2mm

-- Either $a_{ik}=0$, $w=xs_{k}$, $l(w)=l(x)+1$.
Then, $x(\a_i)=\a_j$.
Using (4.6.4) for $x$, and (4.6.3), we get
$$S_w^*(\Fc_j)=q^{a(x,i)}S_{k}^*(F_i^{(x^{-1}(\mu))})=
q^{a(x,i)}F_i^{(\mu)}.$$
Using the identity
$$\Delta_+\cap w^{-1}\Delta_-=s_k(\Delta_+\cap x^{-1}\Delta_-)\cup\{\a_k\}$$
we get also $a(w,i)=a(x,i).$ 
Thus (4.6.4) holds.

\vskip2mm

-- Either $a_{ik}=-1$, $w=xs_is_{k}$, $l(w)=l(x)+2$.
Then, $x(\a_{k})=\a_j$.
Using (4.6.4) for $x$ we get
$$S_w^*(\Fc_j)=q^{a(x,k)}S_{k}^*S_i^*(F_{k}^{(x^{-1}(\mu))}).$$
We are reduced to the $A_2$ case. 
Set $\nu=x^{-1}(\mu)$.
We have  $w^{-1}(\mu)=s_ks_i(\nu)<s_i(\nu)<\nu$.
A direct computation using (4.6.3) gives
$$S_{k}^*S_i^*(F_{k}^{(\nu)})=qF_i^{(s_ks_i(\nu))}.$$
Using the identity
$$\Delta_+\cap w^{-1}\Delta_-=s_ks_i(\Delta_+\cap x^{-1}\Delta_-)
\cup\{\a_k,\a_i+\a_k\}$$
we get also $a(w,i)=a(x,k)+1.$ 
Thus (4.6.4) holds.

\vskip2mm

Setting $w,i,j\to w_0s_\iu,\iu,i$ in (4.6.4) 
and using the formula for $S_{\iu}^*$, we get
$S_{w_0}^*(\Fc_i)=-q^{a(w_0s_\iu,\iu)+2}F_\iu^{(w_0(\mu))}.$
Thus
$$\o^*(\Fc_i)=-q^{a(w_0s_\iu,\iu)+2}\Fc_\iu^\dagger,\leqno(4.6.5)$$
see (4.6.1). 
Moreover we have, see 3.2,
$$\matrix
a(w_0s_\iu,\iu)
&={1\over 2}\sum_{\a\in\Delta_+\setminus\{\a_\iu\}}(\a_\iu,\a)^2\hfill\cr
&={1\over 4}\kappa(\a_\iu,\a_\iu)-2\hfill\cr
&=\cb-2.\hfill
\endmatrix$$

By definition we have 
$q\Fc_i=\Wc_i-\sum_j[a_{ij}]\Vc_j.$
Using $(4.6.1)$, $(4.6.5)$ and the equality $a_{ij}=a_{\iu\ju}$ 
we get the identity 
$$\sum_j[a_{ij}]\o^*(\Vc_j)=\Wc_i^\dagger+q^{\cb}\Wc_\iu^\dagger-
q^{\cb}\sum_j[a_{ij}]\Vc_\ju^\dagger.$$
Claim 1 is proved.

\vskip2mm

From the definition of the operator $S_i$ we get 
$\dagger S_i=S_i\dagger$ for all $i$.
Claim 2 follows immediately.

\vskip2mm

We now prove Claim 3. 
Using Claim 2 it is sufficient to prove that 
$\o(F_{\l\a})\subseteq F_{\l,w_0*\a}$.
Assume that $\a'=s_i*\a$, $\mu'=s_i(\mu)$ as above.
It is sufficient to prove that
$S_i(F_{\l\a}^{(\mu)})\subseteq F_{\l\a'}^{(\mu')}.$
By \cite{10,\S 1.3} the ring of $G_\a$-invariant polynomials
on $m_{\l\a}^{-1}(0)$ is generated by the following two types of functions :

\vskip2mm

$(i)\quad\tr_{V_j}(B_{h_1}B_{h_2}\cdots B_{h_n})$ for any sequence
$h_1,h_2,...,h_n\in H$ such that $j=h'_1$, $h''_1=h'_2$,...,
$h''_{n-1}=h'_n$, $h''_n=j$,

\vskip2mm

$(ii)\quad\varphi(q_jB_{h_1}B_{h_2}\cdots B_{h_n}p_k)$ for any sequence
$h_1,h_2,...,h_n\in H$ such that $j=h'_1$, $h''_1=h'_2$,...,
$h''_{n-1}=h'_n$, $h''_n=k$, and any linear form $\varphi$ on 
$\Hom(W_k,W_j)$.

\vskip2mm

\noindent
We may assume that $m_i<0$.
Fix an element $(x,x')$ in $Z^\mu_i$.
Set $x=(B_h,p_j,q_j)$, $x'=(B'_h,p'_j,q'_j)$.
In particular we have 
$$B_h=B'_h\ \roman{if}\ h',h''\neq i;\quad
B_{h_1}B_{h_2}=B'_{h_1}B'_{h_2}\ \roman{if}\ h''_1=h'_2=i.$$
Thus any function of type $(i)$ coincide on $x$ and $x'$.
We have also
$$q_j=q'_j,\,p_j=p'_j\ \roman{if}\ j\neq i;\quad
q_ip_i=q'_ip'_i;$$
$$q_iB_{h}=q'_iB'_h\ \roman{if}\ h'=i;\quad
B_{h}p_i=B'_hp'_i\ \roman{if}\ h''=i.$$
Thus any function of type $(ii)$ coincide on $x$ and $x'$.
In particular $x\in F_{\l\a}^{(\mu)}$ iff $x'\in F_{\l\a '}^{(\mu ')}.$
We are done.
\quad\qed
\enddemo

\remark{ Remark} 
The dual of the $\tilde G_\l$-bundle $E^{(\mu)}$ on $Q_\l^{(\mu)}$ is 
$(E^*)^{(\mu)}$, 
where $E^*$ is the dual module,
obtained by composing the $\tilde G_\l$-action by the group automorphism
$(z,g_\l,g_\a)\mapsto(z^{-1},{}^tg_\l^{-1},{}^tg_\a^{-1}).$
Note that, in the particular case where $E=V_i$, $W_i$ 
we have $\Vc_i^\dagger=\Vc_i^*$, $\Wc_i^\dagger=\Wc_i^*$.
\endremark
\vskip3mm

\noindent{\bf Convention.} 
Put $1_{\l\a}=\Oc_{F_{\l\a}}$, $1'_{\l\a}=\Oc_{Q_{\l\a}}$.
It is convenient to set also $1_\l=1_{\l 0}$, $1'_\l=1'_{\l 0}$.
To simplify the notations we put $\nu=w_0*0$.

\vskip3mm

\head 5. The involution on the convolution algebra\endhead
\subhead 5.1\endsubhead
Given smooth quasi-projective $G$-varieties $X_1,X_2,X_3,$
consider the projection $p_{ab}\,:\,X_1\times X_2\times X_3\to X_a\times X_b$
for all $1\leq a,b\leq 3$, $a\neq b$. 
Fix closed subvarieties 
$Z_{ab}\subset X_a\times X_b$
such that the restriction of $p_{13}$ to 
$p_{12}^{-1}Z_{12}\cap p_{23}^{-1}Z_{23}$ is proper and maps to $Z_{13}$. 
The convolution product is the map
$$\star\,:\quad\Kb^G(Z_{12})\times\Kb^G(Z_{23})\to\Kb^G(Z_{13}),\quad
(\Ec,\Fc)\mapsto p_{13\,*}\bigl((p_{12}^*\Ec)\otimes(p_{23}^*\Fc)\bigr).$$
If $Z_{12}=Z_{23}=Z_{13}=Z$, the map $\star$ 
endows $\Kb^G(Z)$ with the structure of an $\Rb^G$-algebra.
See \cite{4} for more details.

\subhead 5.2\endsubhead
Let $\DD_{X_a}$ be the Serre-Grothendieck duality operator on $\Kb^G(X_a)$.
Assume that $X_a$ is connected.
Let $\Omega_{X_a}$ be the canonical bundle of $X_a$,
and let $\Oc_{X_a}$ be the structural sheaf. 
We have 
$$\DD_{X_a}(\Ec)=(-1)^{\dim X_a}\Ec^*\otimes\Omega_{X_a}$$
for any $G$-equivariant locally free sheaf $\Ec$ on $X_a$.
Assume that there is a character $q$ of the group $G$ 
such that $\Omega_{X_a}=q^{-\dim X_a}$ for all $a$.  
Consider the operator $D_{Z_{ab}}=q^{d_{ab}}\DD_{Z_{ab}},$
where $d_{ab}=(\dim X_a+\dim X_b)/2.$
Recall that the automorphism $\phi\,:\,X_a\times X_b\to X_b\times X_a$
is the flip.

\proclaim{Lemma}
Fix $x\in\Kb^G(Z_{12})$, $y\in\Kb^G(Z_{23})$. 
\roster
\item $\phi^*(x\star y)=\phi^*(y)\star\phi^*(x),$
$D_{Z_{12}}(x)\star D_{Z_{23}}(y)=D_{Z_{13}}(x\star y),$
$\phi^*D_{Z_{ab}}=D_{Z_{ba}}\phi^*$.

\item If $Z_{12}=Z_{23}=Z_\l$ then
$(\o\times\o)^*(x)\star(\o\times\o)^*(y)=(\o\times\o)^*(x\star y).$

\item If $Z_{12}=Z_\l$, $Z_{23}=Q_\l$ or $F_\l$ then
$(\o\times\o)^*(x)\star\o^*(y)=\o^*(x\star y).$
\endroster
\endproclaim

\noindent
See \cite{14} for more details.

\subhead 5.3\endsubhead
We consider the maps $\g_\l$, $\g'_\l$, $\G_\l$, $\zeta_\l$ 
on $\Kb^{\tilde G_\l}(F_\l)$, $\Kb^{\tilde G_\l}(Q_\l)$,
$\Kb^{\tilde G_\l}(Z_\l)$ such that
$$\g_\l=\Oplus_\a q^{d_{\l\a}/2}\o^*\DD_{F_{\l\a}},\quad
\g'_\l=\Oplus_\a q^{3d_{\l\a}/2}\o^*\DD_{Q_{\l\a}},$$
$$\G_\l=\Oplus_{\a,\a'}(\o\times\o)^*D_{Z_{\l\a\a'}},\quad
\zeta_\l=(\o\times\o)^*\phi^*$$
(see Lemma 4.6.3 for $\g_\l$).
Let 
$$\bar{}\,:\,
\Rb^{\tilde G_\l}\to\Rb^{\tilde G_\l},
\ V\mapsto\bar V$$
be the ring automorphism induced by the group automorphism
$\tilde G_\l\to\tilde G_\l$, $(z,g_\l)\mapsto (z^{-1},g_\l)$.
By Lemma 4.6.1 the operators $\o^*$, $\zeta_\l$ are $\dagger$-semilinear
automorphisms of $\Rb^{\tilde G_\l}$-modules,
and $\g_\l$, $\g'_\l$, $\G_\l$ are $\bar{\ }$-semilinear.
Let $\kappa\,:\,F_\l\hookrightarrow Q_\l$ be the closed embedding.

\proclaim{Lemma} The following identities hold :
\roster
\item $\o^*\DD_{F_\l}=\DD_{F_\l}\o^*$, $\o^*\DD_{Q_\l}=\DD_{Q_\l}\o^*$,
$(\o\times\o)^*\DD_{Z_\l}=\DD_{Z_\l}(\o\times\o)^*$, 

\item $\kappa_*\o^*=\o^*\kappa_*$,
$(\o\times\o)^*\phi^*=\phi^*(\o\times\o)^*$,
$(\kappa_*\times\kappa_*)\phi^*=\phi^*(\kappa_*\times\kappa_*)$,

\item $\g_\l(u\star x)=\G_\l(u)\star\g_\l(x),$
for any $x\in\Kb^{\tilde G_\l}(F_\l)$,
$u\in\Kb^{\tilde G_\l}(Z_\l)$,

\item $\g'_\l(u\star x)=q^{d_{\l\a}-d_{\l\a'}}\G_\l(u)\star\g'_\l(x),$
for any $x\in\Kb^{\tilde G_\l}(Q_{\l\a'}),$
$u\in\Kb^{\tilde G_\l}(Z_{\l\a\a'})$.
\endroster
\endproclaim

\subhead 5.4\endsubhead
Let $\Ab_{\l\a\a'}$ be the quotient of the
$\Rb^{\tilde G_\l}$-module $\Kb^{\tilde G_\l}(Z_{\l\a\a'})$
by its torsion submodule.
We set $\Ab_\l=\Oplus_{\a,\a'}\Ab_{\l\a\a'}$.
Setting $Z_{12}=Z_{23}=Z_\l$ in 5.1, we get
an associative product on the space $\Ab_\l$.
The rings $\Rb^{\CC^\times}$, $\AA$ are identified as in 4.5.
An $\AA$-algebra homomorphism $\Phi_\l\,:\,\Ub\to\Ab_\l$ is given in \cite{20}. 
In this subsection we fix a particular normalization for $\Phi_\l$.
Let $\delta\,:\,Q_\l\hookrightarrow Q_\l\times Q_\l$
be the diagonal embedding, 
and let $p,p'\,:\,Q_\l\times Q_\l\to Q_\l$ 
be the first and the second projection.
Let $f_{i;\a},$ $f^\pm_{i;\a}$, $v_{i;\a}$
be the ranks of $\Fc_{i;\a}$, $\Fc^\pm_{i;\a},$ $\Vc_{i;\a}.$
We have 
$$f_{i;\a}=(\a_i,\l-\a),\quad
f^-_{i;\a}=-v_{i;\a}=-(\o_i,\a).$$
Set 
$t_\a=(\a,2\l-\a)/2+|\a|^2$,
and
$$r^+_{i;\a}=(\l,\a_i)-\bigl(\o_i-\Sum_jn_{ij}\o_j,\a\bigr),\quad
r^-_{i;\a}=-\bigl(\o_i-\sum_j\bar n_{ij}\o_j,\a\bigr).$$
Let $1_{\l\a\a'}\in\Ab_\l$ be the class of the 
structural sheaf of $X_{\l\a\a'}$.
For any $r\in\ZZ$ we put 
$$x_{ir}^+=q^{(1-\cb)r}\sum_{\a'=\a+\a_i}
(-1)^{r^+_{i;\a'}}
(q^{-1}\Wedge_{\Vc}^{-1}\boxtimes\Wedge_{\Vc})^{r+f^-_{i;\a'}}
\otimes {p'}^*\Wedge_{\Fc^+_{i}}^{-1}\otimes\Wedge_\Wc^{t_{\a'}-t_\a}
\otimes 1_{\l\a\a'},$$
$$x_{ir}^-=q^{(1-\cb)r}\sum_{\a'=\a-\a_i}
(-1)^{r^-_{i;\a'}}
(q^{-1}\Wedge_{\Vc}\boxtimes\Wedge_{\Vc}^{-1})^{r+f^+_{i;\a'}}
\otimes {p'}^*\Wedge_{\Fc^-_{i}}^{-1}\otimes\Wedge_\Wc^{t_{\a'}-t_\a}
\otimes 1_{\l\a\a'}.$$
Let also $k_i^\pm(z)$ be the expansion at $z=\infty$ or $0$ of 
$$\delta_*\sum_\a q^{f_{i;\a}}
\bigl(\Sum_{r\geq 0}(-q^{-\cb}/z)^r\wedge^r{\Fc_{i;\a}}\bigr)\otimes
\bigl(\Sum_{r\geq 0}(-q^{2-\cb}/z)^r\wedge^r{\Fc_{i;\a}}\bigr)^{-1}.$$
The map $\Phi_\l$ takes 
$\xb^\pm_{ir}$ to $x^\pm_{ir},$ 
and $\kb^\pm_{ir}$ to $k^\pm_{ir}.$ 
For a future use, let us mention the following.

\proclaim{Lemma} For any $n>0$ we have
$$(x_{i0}^+)^{(n)}=\pm\sum_{\a'=\a+n\a_i}
(q^{-n}\Wedge_{\Vc}^{-1}\boxtimes\Wedge_{\Vc})^{f^-_{i;\a'}}
\otimes {p'}^*\Wedge_{\Fc^+_{i}}^{-n}\otimes\Wedge_\Wc^{t_{\a'}-t_\a}
\otimes 1_{\l\a\a'},$$
and similarly for $(x_{i0}^-)^{(n)}$.
\endproclaim

\demo{Proof}
By the same argument as in \cite{20, \S 11.1-3} 
it is enough to check this relation for type $A_1$.
In this case, using the faithful representation introduced in
\cite{22} the formula follows from a direct computation :
the formula for the action of
the operator $(x_{i0}^+)^{(n)}$
may be found in \cite{20, Lemma 12.1.1}, 
its relation with the identity above is proved
as in \cite{22}.
\quad\qed
\enddemo
\vskip3mm

\noindent{\bf Convention.}
Hereafter we may omit the maps $\delta_*,p^*,{p'}^*,\otimes$, 
hoping that it makes no confusion.

\subhead 5.5\endsubhead
Let $H_\l\subseteq G_\l$ be any closed subgroup.
Set $\tilde H_\l=H_\l\times\CC^\times$.
For simplicity we set
$\Wb_{H_\l,\a}=\Kb^{\tilde H_\l}(F_{\l\a})$, 
$\Wb'_{H_\l,\a}=\Kb^{\tilde H_\l}(Q_{\l\a})$,
$\Wb_{H_\l}=\oplus_\a\Wb_{H_\l,\a}$,
$\Wb'_{H_\l}=\oplus_\a\Wb'_{H_\l,\a}$.
Taking $Z_{12}=Z_\l$, $X_3=\{point\}$, and
$Z_{13}=Z_{23}=Q_\l$ or $F_\l$ in 5.1 we get a left
$\Ub$-action on the $\Rb^{\tilde H_\l}$-modules
$\Wb_{H_\l}$, $\Wb'_{H_\l}$ such that
$$(u,x)\mapsto u\cdot x=\Phi_\l(u)\star x.$$
Taking $Z_{23}=Z_\l$, $X_1=\{point\}$, and
$Z_{12}=Z_{13}=Q_\l$ or $F_\l$ we get a right 
$\Ub$-action on $\Wb_{H_\l}$, $\Wb'_{H_\l}$ such that
$x\cdot u=x\star\Phi_\l(u)=\phi^*\Phi_\l(u)\star x,$
see Lemma 5.2.
We fix a maximal torus $T_\l\subset G_\l$.

\proclaim{Lemma}
\roster
\item The $\Rb^{\tilde H_\l}$-modules 
$\Wb_{H_\l}$, $\Wb'_{H_\l}$
are free of finite type, and we have
$\Wb_{H_\l}=\Wb_{G_\l}\otimes_{\Rb_\l}\Rb^{\tilde H_\l},$
$\Wb'_{H_\l}=\Wb'_{G_\l}\otimes_{\Rb_\l}\Rb^{\tilde H_\l}.$
Moreover, there is a canonical action of the Weyl group of $G_\l$ on
$\Wb_{T_\l}$, $\Wb'_{T_\l}$ such that the forgetful map
identifies $\Wb_{G_\l}$, $\Wb'_{G_\l}$ with the subspaces of invariant
elements in $\Wb_{T_\l}$, $\Wb'_{T_\l}$.

\item We have $\Wb_{H_\l}=\Ub\cdot(\Rb^{H_\l}\otimes 1_\l).$

\item We have $\Wb_{G_\l}=\Ub\cdot 1_\l.$
\endroster
\endproclaim

\demo{Proof}
Claim 1 is proved in \cite{20, Theorem 7.3.5}.
See also \cite{4, Chapter 5}.
Claim 2 is proved as in \cite{20, Proposition 12.3.2}.
Let us prove that, if $H_\l=G_\l$, then 
$\Rb^{G_\l}\otimes 1_\l\in\Ub\cdot 1_\l$. 
By construction $\kb_i^\pm(z)\cdot 1_\l$ is the expansion of
$$q^{f_{i;0}}
\bigl(\Sum_{r\geq 0}(-1/qz)^r\wedge^rW_i\bigr)\otimes
\bigl(\Sum_{r\geq 0}(-q/z)^r\wedge^rW_i\bigr)^{-1}$$
in $\Rb^{\tilde G_\l}[[z^{-1}]]$. 
Thus, the elements $\pb_{i,s}\in\Rb^{\GL_{\ell_i}}$ 
are the elementary symmetric polynomials or zero. 
Claim 3 follows.
\quad\qed
\enddemo
\vskip3mm

\noindent{\bf Convention.}
Although most of our constructions are meaningful
for any closed subgroup $H_\l\subseteq G_\l$, hereafter $H_\l$ 
will be either $G_\l$ or $T_\l$.
For simplicity, if $H_\l=G_\l$ we put $\Wb_\l=\Wb_{G_\l}$, $\Wb'_\l=\Wb'_{G_\l}$,
$\Rb_\l=\Rb^{\tilde G_\l}$, $\Xb_\l=\Xb^{G_\l}$. 

\subhead 5.6\endsubhead
The same construction as in 5.5 yields also
a left and a right action of $\Ab_\l$ on
$\Wb_{H_\l}$, $\Wb'_{H_\l}.$
Since $\Wb_{H_\l}$, $\Wb'_{H_\l}$
are integrable left $\Ub$-modules,
they admit a left $\Rb^{\tilde H_\l}$-linear action of the group $B_W$.
We normalize this action in such a way that the element $T_i\in B_W$ 
acts as the operator $T''_{i,1}$ in \cite{13, \S 5.2.1}.
Similarly let $\check T_w$ be the left action of the element
$T_w\in B_W$ associated to the right $\Ub$-action.

\proclaim{Proposition}
\roster
\item There is a unique action of the group $B_W$ on $\Ab_\l$
by $\Rb^{\tilde H_\l}$-algebra automorphisms such that
$$T_w(x\star y)=T_w(x)\star T_w(y),\quad
\forall x\in\Ab_\l, y\in\Wb_{H_\l}\ \text{or}\ \Wb'_{H_\l},\ w\in W.$$

\item We have 
$T_w\Phi_\l=\Phi_\l T_w,$ and
$\check T_w(y\star x)=\check T_w(y)\star\sigma(T_w)(x)$
for any $x,y$ as above. 

\item There is an invertible element $r_\l\in\Rb^{\tilde H_\l}$
such that 
$T_{w_0}(1'_{\l\nu})=r_\l\otimes 1'_{\l}$,
$T_{w_0}(1_{\l\nu})=r_\l\otimes 1_{\l}$.

\item There is an invertible element $s_\l\in\Rb^{\tilde H_\l}$
such that 
$T_{w_0}(1'_\l)=s_\l\otimes 1'_{\l\nu}$,
$T_{w_0}(1_\l)=s_\l\otimes 1_{\l\nu}$.

\item There is an invertible element $\vartheta\in\AA$ such that
$r_\l\otimes 1_{\l\nu}=\vartheta\otimes\Wedge_\Wc^{t_\nu}\otimes
\Otimes_i\Wedge^{-v_{i;\nu}}_{\Vc_{i;\nu}}.$
Moreover $r_\l s_\l=(-q)^{(\rho,\nu)}$ and
$\vartheta=\pm q^{|\nu|^2/2+\cb(\l,\l)/2}$.
\endroster
\endproclaim

\demo{Proof}
Let us prove Claim 1.
We fix elements $i\in I$, $\l\in P^+$.
Let $\Ub'_i\subset\Ub'$ be the subalgebra
generated by $\eb_i$, $\fb_i$, $\kb_i^{\pm 1}$.
The proof uses an element $\tau'_i$ introduced in \cite{17}.
The element $\tau'_i$ belongs to a ring completion of $\Ub'_i$, 
is invertible, and satisfies the following identity 
$$\tau'_i\cdot x\cdot{\tau_i'}^{-1}=T'_{i,1}(x),\quad\forall x\in\Ub.$$
Let us recall the construction of $\tau'_i$, following \cite{5}.

For any $\QQ(q)$-vector space $V$ we set $V^*=\Hom_{\QQ(q)}(V,\QQ(q))$.
Let $\Rb'_i\subset{\Ub'_i}^*$ be the 
$\QQ(q)$-space spanned by the matrix elements
of the finite dimensional $\Ub'_i$-modules.
It is a Hopf algebra.
Let $\sum x_0\otimes x_1$ denote the image of the element $x\in\Rb'_i$
by the coproduct.
The space ${\Rb'_i}^*$ is a ring such that 
$(f\cdot g)(x)=\sum f(x_0)\,g(x_1)$ for all $x\in\Rb'_i$.
The canonical map $\Ub'_i\to{\Rb'_i}^*$ is a ring homomorphism.
An integrable $\Ub'$-module $V$ restricts to an integrable $\Ub'_i$-module
via the canonical embedding $\Ub'_i\subset\Ub'$.
It is also a $\Rb'_i$-comodule for the co-action
$V\to\Rb'_i\otimes V$, $v\mapsto\sum v_1\otimes v_2$ such that 
$x\cdot v=\sum v_1(x)\,v_2$ for all $x\in\Ub'_i$.
Therefore the ring ${\Rb'_i}^*$ acts on $V$ by
$f\cdot v=\sum f(v_1)\,v_2.$
This action restricts to the original $\Ub'_i$-action via the 
map $\Ub'_i\to{\Rb'_i}^*$.
Let $\tau_i\in{\Rb'_i}^*$ be the element
denoted by $t$ in \cite{5, \S 1.6}.
It is invertible.
Let $t_i$ be the operator on $V$ 
taking $v$ to $\tau_i\cdot v$.
It is invertible and the inverse takes $v$ to $\tau_i^{-1}\cdot v$.
We have
$$t_i(v)=T''_{i,-1}(v),\qquad
t_ixt_i^{-1}(v)=T''_{i,-1}(x)\cdot v,$$ 
for all $x\in\Ub'$ and $v\in V$, by \cite{5, \S 2}.
Therefore
$$\tau_i\cdot x\cdot{\tau_i}^{-1}=T''_{i,-1}(x),\quad
\forall x\in\Ub',$$
where we take the product in the ring ${\Rb'_i}^*$ in the left hand side.

For any $n\in\NN$ let
$\Lambda_i(n)$ be the simple 
$\Ub'_i$-module with highest weight $n\o_i$.
Let $\Rb'_{in}\subset\Rb'_i$ be the subspace
spanned by matrix elements 
of the module $\Oplus_{n'\leq n}\Lambda_i(n')$.
It is a subcoalgebra.
Let $\Ib'_{in}\subset\Ub'_i$ be the annihilator
of $\Oplus_{n'\leq n}\Lambda_i(n')$.
It is a two-sided ideal.
The canonical map 
$\Ub'_i\to{\Rb'_i}^*$
factorizes through an isomorphism 
$\Ub'_i/\Ib'_{in}\simto(\Rb'_{in})^*.$
Since
$\Rb'_i={\lim\limits_\lra}_n\Rb'_{in},$
we get a $\QQ(q)$-algebra isomorphism
$$
{\lim\limits_\lla}_n(\Ub'_i/\Ib'_{in})
\simto
{\Rb'_i}^*
.$$
For each $n$ we choose $\tau_{in}\in\Ub'_i$ such that
$$\tau_{in}-\tau_i\in{\lim\limits_\lla}_{n'\geq n}
(\Ib'_{in}/\Ib'_{in'}).$$

Since $\Ub'_i$ embeds in $\Ub'$,
the space $\Ab_\l\otimes_\AA\QQ(q)$ 
is a $\Ub'_i\otimes_\AA\Rb^{\tilde H_\l}$-bimodule
of finite type over $\Rb^{\tilde H_\l}\otimes_\AA\QQ(q)$.
In particular there is an integer $n$ such that 
the ideal $\Ib'_{in}$ acts trivially on $\Ab_\l\otimes_\AA\QQ(q)$.
Fix such an integer $n$.
Then the operator $T''_{i,-1}$ acts on $\Ab_\l\otimes_\AA\QQ(q)$
via the conjugation by the element 
$\Phi_\l(\tau_{in})\in\Ab_\l\otimes_\AA\QQ(q).$ 
Moreover, the formulas in \cite{5, \S 2} imply that
the left and right product by
$\Phi_\l(\tau_{in})$, $\Phi_\l({\tau_{in}}^{-1})$ preserves the subspace
$\Ab_\l\subset\Ab_\l\otimes_\AA\QQ(q)$.

Recall that $T_i$ acts as Lusztig's operator $T''_{i,1}$.
The element $\tau'_i$ yielding the
action of $T_{i,1}''$ on $\Ab_\l$ 
can be constructed as $\tau_i$,
using the identity 
$$T_{i,1}''(x)=(-q)^{(\a,\a_i)}(T''_{i,-1})^{-1}(x)$$
for any $x\in\Ub_\a$, see \cite{13, \S 37.2.4}.

Recall that $N_{\l}$ is a cone over the point 0 equal
to the class of the trivial representation, 
and that the fixed points subset $(N_\l)^{\CC^\times}$ is reduced to $\{0\}$.
Hence $A_\l\otimes_\AA\QQ(q)$ coincides with the tensor product
$(\Wb_{G_\l}\otimes_{\Rb^{\tilde G_\l}}\Wb_{G_\l})\otimes_\AA\QQ(q)$ 
by the Kunneth isomorphism, see \cite{4, 5.6} and \cite{20,\S 7}.
Since $A_\l$ is torsion free over $\AA$, it embeds in $A_\l\otimes_\AA\QQ(q)$.
Then, a standard argument implies that
$\Wb_{G_\l}$ is a faithful $\Ab_\l$-module,
see \cite{4, \S 5} for more details.
For the same reason $\Wb'_{H_\l}$ is also a faithful $\Ab_\l$-module.
The unicity in Claim 1 follows.

\vskip2mm

Claim 2 is obvious from the previous construction.
Claims 3,4  are obvious either since 
$T_{w_0}$ is an invertible $\Rb^{\tilde H_\l}$-linear
homomorphism from $\Wb'_{H_\l,0}$ to $\Wb'_{H_\l,\nu}$,
and both $\Rb^{\tilde H_\l}$-modules are free of rank one,
generated by $1'_\l,1'_{\l\nu}$ respectively.

\vskip2mm

Let us prove Claim 5.
Part two of the claim follows from 
\cite{11, 5.4.$(a)$ and Corollary 5.9}
(note that the formula for $s'$ in 
\cite{11, 5.4.$(a)$} should be replaced by 
$s'=\sum_{i,j}a_{ij}t_it_j/2-\sum_i t_i(1+d_i)$).
Let us prove Part one.
Consider an element $w\in W$ such that $l(w)>0$.
Fix $i\in I$ such that $l(s_iw)=l(w)-1$.
Put $w'=s_iw$.
Set $\a'=w'*\nu,$ $\a=w*\nu,$ $n=(ww_0(\l),\a_{i}).$
Thus $n>0$ and $\a=\a'-n\a_{i}$.
We have 
$$T_{i}(1'_{\l\a'})=\eb_{i}^{(n)}(1'_{\l\a'}),$$
see \cite{11, Lemma 5.6}.
Let $r_w\in\Rb^{\tilde H_\l}$ be the unique element such that
$T_w(1'_{\l\nu})=r_w1'_{\l\a}.$
The varieties $Q_{\l\a}$, $Q_{\l\a'}$ are reduced to a point.
Thus, using Lemma 5.4 we get
$$T_w(1'_{\l\nu})=\pm
\bigl(q^{-n}\Wedge^{-1}_{\Vc_{\a}}\Wedge_{\Vc_{\a'}}\bigr)
^{f^-_{i;\a'}}\Wedge^{-n}_{\Fc^{+}_{i;\a'}}\,
\Wedge_\Wc^{t_{\a'}-t_\a}
(r_{w'}\,1'_{\l\a}).$$
The classes of the $\tilde H_\l$-equivariant sheaves 
$\Vc_{i;\a},\Vc_{i;\a'},\Fc^+_{i;\a'}$ 
are identified with elements
of $\Rb^{\tilde H_\l}$ in the obvious way.
Let first assume that we have 
$$\Wedge_{\Fc^{+}_{i;\a'}}=
\Wedge_{\Vc_{i;\a}}\in\Rb^{\tilde H_\l}.\leqno(5.6.1)$$
Then
$$r_w=\pm q^s\Wedge_{\Vc_{i;\a'}}^{-v_{i;\a'}}
\Wedge^{v_{i;\a}}_{\Vc_{i;\a}}\Wedge_\Wc^{t_{\a'}-t_\a} r_{w'},
\ \roman{where}\ 
s=(\a'-\a,\o_i)(\a',\o_i)
.$$
By induction on $l(w)$ we get
$$r_{w}\in\pm q^\ZZ\Wedge_\Wc^{t_\nu-t_\a}
\prod_j\Wedge_{\Vc_{j;\nu}}^{-v_{j;\nu}}
\Wedge_{\Vc_{j;\a}}^{v_{j;\a}}.\leqno(5.6.2)$$
Setting $w=w_0$ in (5.6.2) we get
$r_\l=\vartheta\Wedge_\Wc^{t_\nu}\Prod_i\Wedge^{-v_{i;\nu}}_{\Vc_{i;\nu}},$
with $\vartheta\in\pm q^\ZZ.$
A direct computation gives
$\vartheta\in\pm q^t$, where
$$\matrix
t
&=|\nu|^2/2+\sum_{\a\in\Delta_+}(\l,\a)^2/2\hfill\cr
&=|\nu|^2/2+\cb(\l,\l)/2,\hfill
\endmatrix$$
see \cite{6, Exercice 6.2}.

\vskip2mm

Finally, we prove (5.6.1).
We have an isomorphism of $\tilde G_\l$-varieties
$Q_{\l\a}^{(\mu)}\to Q_{\l\a'}^{(\mu')}$, 
see the proof of Lemma 4.6 and the notations therein.
This isomorphism takes $V_i^{(\mu)}$ to $(F_i^+)^{(\mu')}$.
Assume that $\mu$ is regular dominant, so that
$Q_{\l\a}^{(\mu)}=Q_{\l\a}$ and
$V_i^{(\mu)}=\Vc_{i;\a}$.
Since the $\tilde G_\l$-variety $Q_{\l\a'}$ is reduced to a point,
it is canonicaly isomorphic to $Q_{\l\a'}^{(\mu')}$, 
and the isomorphism takes $\Fc^{+}_{i;\a'}$ to $(F_i^+)^{(\mu')}$.
\quad\qed
\enddemo
\subhead 5.7\endsubhead
For each $i,\a$ we consider the elements 
$$
g_{i;\a}=-1+(\cb-1)f^+_{\iu;w_0*\a}+f^+_{i;\a}-\cb f^-_{\iu;w_0*\a}\in\ZZ,
$$
$$
h_{i;\a}=r^+_{\iu;w_0*\a}+d_{\l,w_0*\a-\a_\iu,w_0*\a}+r^-_{i;\a}
\in\ZZ/2\ZZ.
$$

\noindent{\bf Convention.}
The elements $r_{i;\a}^\pm$ depend on the choice of the orientation $\O$.
Hereafter we assume that $n_{ij}=n_{\iu\ju}$ for all $i,j$ if $\cb$ is even,
and $n_{ij}=\bar n_{\iu\ju}$ for all $i,j$ if $\cb$ is odd 
(i.e. if $\gen$ is of type $A_{2n}$).
The existence of such an orientation is checked case-by-case.

\vskip3mm

\noindent
Using the convention above, a direct computation gives
$$g_{i;\a+\a_j}+g_{j;\a}=g_{j;\a+\a_i}+g_{i;\a},\quad
h_{i;\a+\a_j}+h_{j;\a}=h_{j;\a+\a_i}+h_{i;\a}$$
for all $i,j$.
Thus there are unique quadratic maps 
$x\,:\,Q\to\ZZ$, $y\,:\,Q\to\ZZ/2\ZZ$ 
such that
$$\matrix
x(\nu)=(\cb-1)|\nu|^2-\cb(\l,\l),
\quad\hfill&
x(\a+\a_i)-x(\a)=\cb v_{i;\nu}-g_{i;\a},\hfill\cr\cr
y(\nu)=0,
\quad\hfill&
y(\a+\a_i)-y(\a)=h_{i;\a},\hfill
\endmatrix
\leqno(5.7.1)$$
for all $\a$, $i$. Put $\xi(\a)=q^{x(\a)},$ $\varrho(\a)=(-1)^{y(\a)},$
$t_{\a\a'}=t_{w_0*\a'}-t_{\a'}-t_{w_0*\a}+t_\a.$
Consider the element
$$c_{\l\a}=\varrho(\a)
\xi(\a)q^{-\cb |\a|^2}
\Wedge^{t_{0\a}}_\Wc\otimes
\Otimes_i\left(
\Wedge_{\Vc_{i;\a}}^{2v_{i;\nu}}\otimes
\Wedge_{\Vc_{i;\a}^*\otimes\o^*\Vc_{\iu;w_0*\a}^*}\right)
\otimes r_\l
\in\Wb'_{H_\l,\a}.$$
Put $c_\l=\sum_\a c_{\l\a}$.

\proclaim{Lemma} 
\roster
\item If $\a'=\a-\a_j$ then the restriction of
$c_{\l}^{-1}\boxtimes c_{\l}$ to $X_{\l\a\a'}$ is
$$(-1)^{h_{j;\a'}}q^{g_{j;\a'}}
(\Wedge_{\Vc}\boxtimes\Wedge_{\Vc}^{-1})
^{-f^+_{j;\a'}-f^+_{\ju;w_0*\a'}}
\Wedge^{t_{\a\a'}}_\Wc
\otimes{p'}^*\Wedge_{\Fc_j^-+q^\cb\o^*\Fc_\ju^{-*}}
\otimes 1_{\l\a\a'}.$$

\item We have $\o^*(c_{\l\a})=c_{\l,w_0*\a}$, 
$c_{\l\nu}=r_\l^{-1}\otimes 1'_{\l\nu}$.
\endroster
\endproclaim

\demo{Proof}
Fix $\a,\a'$ such that $\a'=\a-\a_j$. 
For any $i$ we consider the following elements
in $\Wb'_{H_\l}$ :
$$\Uc_{i}=q^\cb\o^*(\Vc_{\iu}^*)+\Vc_{i},\quad
c'_{\l}=\Otimes_i\Wedge_{q^{\cb}\Vc_{i}^*\otimes\o^*\Vc_{\iu}^*}.$$
The rank of $\Uc_{i;\a}$ is $v_{i;\a}+v_{\iu;w_0*\a}=v_{\iu;\nu}$.
From Lemma 4.6.1 we have
$$\sum_j[a_{ij}]\Uc_{j;\a}=q^\cb\Wc_{\iu;\a}+\Wc_{i;\a}.$$
The quantum Cartan matrix 
(i.e. the $I\times I$-matrix whose $(i,j)$-th entry is $[a_{ij}]$)
is invertible over $\QQ(q)$.
Thus, for any $\a,\a',i$ we get 
$$(1_{{\l\a}}\boxtimes\Uc_{i;\a'})|_{X_{\l\a\a'}}=
(\Uc_{i;\a}\boxtimes 1_{{\l\a'}})|_{X_{\l\a\a'}}
\in\Kb^{\tilde H_\l}(X_{\l\a\a'}).\leqno(5.7.2)$$
We have $c'_{\l} =\Otimes_i\Wedge_{\Vc_{i}^*\otimes\Uc_{i}}$
and $\Fc_i^-=-\Vc_i$.
Thus, using (5.7.2) we get
$$({c'}_{\l}^{-1}\boxtimes c'_{\l})|_{X_{\l\a\a'}}=
(\Wedge_{\Vc_\a}\boxtimes\Wedge_{\Vc_{\a'}}^{-1})
^{v_{j;\nu}}|_{X_{\l\a\a'}}
\otimes{p'}^*\Wedge_{\Fc_j^-+q^\cb\o^*\Fc_\ju^{-*}}.$$
Note that
$$c_{\l\a}=\xi(\a)\varrho(\a)r_\l\otimes c'_{\l\a}\otimes\Otimes_i
\Wedge_{q^{-\cb/2}\Vc_{i;\a}}^{2v_{i;\nu}}\otimes\Wedge_\Wc^{t_{0\a}}.$$
Thus Claim 1 follows, using (5.7.1) and the identity
$v_{j;\nu}=f^+_{j;\a'}+f^+_{\uj;w_0*\a'}.$
Let us prove Claim 2.
By definition of $c'_\l$ we have $\o^*(c'_{\l\a})=c'_{\l,w_0*\a}$ 
and $c'_{\l 0}=1'_\l$.
Thus, using Proposition 5.6.4 we get
$$c_{\l\nu}= 
\vartheta^{-2}r_\l\otimes\Wedge_\Wc^{-2t_\nu}\otimes\Otimes_i\Wedge_{\Vc_{i;\nu}}^{2v_{i;\nu}}
=r_\l^{-1}\otimes 1'_{\l\nu}
.$$
Since $Q_{\l\nu}$ is a point,
we identify the equivariant sheaf
$\Wedge_{\Vc_{i;\nu}}$
with an element in $\Rb^{\tilde H_\l}$ in the obvious way.
Using (5.7.2) we get
$$\matrix
\Otimes_i
\Wedge_{q^{-\cb/2}\o^*(\Vc_{i;\a})+q^{\cb/2}\Vc^*_{\iu;w_0*\a}}
^{2v_{i;\nu}}
&=
\left(\Otimes_i
\Wedge_{q^{-\cb/2}\o^*(\Vc_{i;0})+q^{\cb/2}\Vc^*_{\iu;\nu}}
^{2v_{\iu;\nu}}
\right)
\otimes 1'_{\l,w_0*\a}\hfill\cr\cr
&=
\vartheta^{-2}q^{|\nu|^2\cb}r_\l^2\otimes\Wedge_\Wc^{-2t_\nu}
\otimes 1'_{\l,w_0*\a}.\hfill
\endmatrix$$
Thus, $\o^*(r_\l)=\vartheta^2q^{-|\nu|^2\cb}r_\l^{-1}$ and
$$\matrix
\o^*\left(
r_\l\otimes\Wedge_\Wc^{-t_\nu}\otimes\Otimes_i\Wedge_{q^{-\cb/2}\Vc_{i;\a}}
^{2v_{i;\nu}}\right)
&=\vartheta^2q^{-|\nu|^2\cb}r_\l^{-1}\otimes\Wedge_\Wc^{t_\nu}\otimes\Otimes_i
\Wedge_{q^{-\cb/2}\o^*\Vc_{i;\a}}^{2v_{i;\nu}}
\hfill\cr\cr
&=r_\l\otimes\Wedge_\Wc^{-t_\nu}\otimes
\Otimes_i\Wedge_{q^{-\cb/2}\Vc_{\iu;w_0*\a}}
^{2v_{i;\nu}}
.\hfill
\endmatrix
$$
A direct computation (see the Appendix) shows that  
$\xi(\a)=\xi(w_0*\a)$,
$\varrho(\a)=\varrho(w_0*\a)$ for all $\a$.
We are done.
\quad\qed
\enddemo

Let $C_\l$ be the $\Rb_\l$-linear automorphism
of $\Ab_\l$ such that 
$$C_\l(x)=x\otimes(c_\l\boxtimes c_\l^{-1}).$$

\proclaim{Proposition}
\roster 
\item The map $C_\l\G_\l$ is an algebra involution of $\Ab_\l$ such that
$C_\l\G_\l(q)=q^{-1}$, $C_\l\G_\l(x_{ir}^\pm)=q^{r\cb}x_{\iu r}^\mp$.

\item The map $C_\l\zeta_\l$ is an algebra anti-involution of
 $\Ab_\l$ such that 
$C_\l\zeta_\l(q)=q$, $C_\l\zeta_\l(x_{ir}^\pm)=q^{-r\cb}x_{\iu,-r}^\pm$.
\endroster
\endproclaim

\demo{Proof}
The variety $X_{\l\a\a'}$ is smooth of dimension 
$d_{\l\a\a'}:=(d_{\l\a}+d_{\l\a'})/2$.
Let $\O_{\l\a\a'}$ be its canonical bundle.
If $\a'=\a+\a_i$, using (4.5.1), Lemma 4.6.1 and Remark 4.6 we get
$$\O_{\l\a\a'}=
q^{f_{i;\a'}-d_{\l\a'}}
(q^{-1}\Wedge_{\Vc_\a}^{-1}\boxtimes\Wedge_{\Vc_{\a'}})^{f_{i;\a'}}
\otimes {p'}^*\Wedge^{-1}_{\Fc_{i;\a'}},\leqno(5.7.3)$$
$$d_{\l\a\a'}-d_{\l\a'}+f_{i;\a'}=-1,\leqno(5.7.4)$$
$$(\o\times\o)^*\bigl(\Wedge_{\Vc_{i;\a}}\boxtimes\Wedge_{\Vc_{i;\a'}}^{-1}\bigr)=
q^{\cb}\Wedge_{\Vc_{\iu;w_0*\a}}\boxtimes\Wedge_{\Vc_{\iu;w_0*\a'}}^{-1}.
\leqno(5.7.5)$$
Using (5.7.3-4) we get
$$\matrix
D_{Z_\l}(x_{ir}^+)=\hfill\cr
={\ds\sum_{\a'=\a+\a_i}
(-1)^{r^+_{i;\a'}+d_{\l\a\a'}}
q^{r\cb-r+d_{\l\a\a'}}
(q^{-1}\Wedge_{\Vc}^{-1}\boxtimes\Wedge_{\Vc})^{-r-f^-_{i;\a'}}
\otimes{p'}^*\Wedge_{\Fc^+_{i}}\otimes\Wedge_\Wc^{t_\a-t_{\a'}}
\otimes\O_{\l\a\a'}}\hfill\cr\cr
={\ds\sum_{\a'=\a+\a_i}
(-1)^{r^+_{i;\a'}+d_{\l\a\a'}}
q^{r\cb-r-1}
(q^{-1}\Wedge_{\Vc}^{-1}\boxtimes\Wedge_{\Vc})^{-r+f^+_{i;\a'}}
\otimes{p'}^*\Wedge_{\Fc^-_{i}}^{-1}\otimes\Wedge_\Wc^{t_\a-t_{\a'}}
\otimes 1_{\l\a\a'}
.}\hfill\cr\cr
\endmatrix$$
Thus, using (5.7.5) we get
$$\matrix
\G_\l(x_{ir}^+)=\hfill\cr
={\ds
\sum_{\a'=\a-\a_\iu}
(-1)^{e_{\iu;\a'}}
q^{-1-r+\cb f_{i;w_0*\a'}^+}
(q\Wedge_{\Vc}\boxtimes\Wedge_{\Vc}^{-1})^{r-f^+_{i;w_0*\a'}}
\otimes{p'}^*\o^*\Wedge_{\Fc^-_{i}}^{-1}\otimes}\hfill\cr
{\qquad\ds\otimes\Wedge_\Wc^{t_{w_0*\a'}-t_{w_0*\a}}\otimes 1_{\l\a\a'}
}\hfill\cr\cr
={\ds 
x_{\iu r}^-\otimes\sum_{\a,\a'}
(-1)^{h_{\iu;\a'}}
q^{r\cb-1+(\cb-2)f_{i;w_0*\a'}^+}
(q^{-1}\Wedge_{\Vc_\a}\boxtimes
\Wedge_{\Vc_{\a'}}^{-1})^{-f^+_{\iu;\a'}-f^+_{i;w_0*\a'}}\otimes}\hfill\cr
{\ds\qquad\otimes\Wedge_\Wc^{t_{\a\a'}}\otimes{p'}^*
\Wedge_{\Fc^-_{\iu}+\o^*\Fc_{i}^{-*}}}\hfill\cr
={\ds 
q^{r\cb}x_{\iu r}^-\otimes\sum_{\a,\a'}
(-1)^{h_{\iu;\a'}}q^{g_{\iu;\a'}}
(\Wedge_{\Vc_\a}\boxtimes
\Wedge_{\Vc_{\a'}}^{-1})^{-f^+_{\iu;\a'}-f^+_{i;w_0*\a'}}\otimes}\hfill\cr
{\ds\qquad\otimes\Wedge_\Wc^{t_{\a\a'}}
\otimes{p'}^*\Wedge_{\Fc^-_\iu+q^\cb\o^*\Fc_i^{-*}},}\hfill\cr
\endmatrix$$
where 
$e_{\iu;\a'}=r^+_{i;w_0*\a'}+d_{\l,w_0*\a,w_0*\a'},$
and $g_{i;\a}$, $h_{i;\a}$ are as at the beginning of  5.7.
Using Lemma 5.7.1 we get
$$C_\l\G_\l(x_{ir}^+)=q^{r\cb}x_{\iu r}^-.$$
Using Lemma 5.2.1 and Lemma 5.3.1-2 we get 
$$\zeta_\l=\G_\l\phi^*D_{Z_\l}=\phi^*D_{Z_\l}\G_\l,\quad
C_\l\phi^*D_{Z_\l}=\phi^*D_{Z_\l}C_\l.$$
Thus, $(C_\l\G_\l)^2=(C_\l\zeta_\l)^2$.
Using Lemma 5.7.2 we get
$$(C_\l\G_\l)^2=(C_\l\zeta_\l)^2=\Id.\leqno(5.7.6)$$
Thus, $C_\l\G_\l(x_{ir}^-)=q^{r\cb}x_{\iu r}^+$ either.
Recall that $\phi^*D_{Z_\l}\Phi_\l=\Phi_\l\tau$, see \cite{24, Lemma 6.5}.
Then Claim 2 follows from Proposition 3.2.2 and (5.7.6).
\quad\qed
\enddemo
\subhead 5.8\endsubhead
Let $A_\l,B_\l\,:\,\Ab_\l\to\Ab_\l$ be the $\Rb_\l$-algebra 
automorphisms such that
$$A_\l(x)=(-q)^{(\rho,\a-\a')}x,\quad
B_\l(x)=(-q)^{(\rho,\a-\a')}q^{(\a'-\a,2\l-\a'-\a)/2}x,$$
for any element $x\in\Ab_{\l\a\a'}$.
Then
$$\Phi_\l A=A_\l\Phi_\l,\quad\Phi_\l B=B_\l\Phi_\l.\leqno(5.8.1)$$
We consider the automorphisms $\b_{Z_\l}$, $\psi_{Z_\l}$
of the ring $\Ab_\l$ such that
$$\b_{Z_\l}=T_{w_0}B_\l C_\l\G_\l,\quad
\psi_{Z_\l}=T_{w_0}A_\l C_\l\zeta_\l.\leqno(5.8.2)$$

\proclaim{Corollary} 
\roster
\item The map $\b_{Z_\l}$ is $\bar{\ }$-semilinear,
the map $\psi_{Z_\l}$ is $\dagger$-semilinear.
Moreover, we have $\b^2_{Z_\l}=\psi_{Z_\l}^2=\Id$.
 
\item For any $u\in\Ub$ we have 
$\Phi_\l(\bar u)=\b_{Z_\l}\Phi_\l (u)$,
$\Phi_\l\psi(u)=\psi_{Z_\l}\Phi_\l(u)$.
\endroster
\endproclaim

\vskip3mm

\demo{Proof}
From Proposition 5.7.2 the map $C_\l\zeta_\l$
is an antihomomorphism of $\Ab_\l$ such that
$q\mapsto q$, $x^\pm_{i0}\mapsto x^\pm_{\iu 0}$ for all $i$.
Thus, using \cite{13, \S 37} we get
$$A_\l C_\l\zeta_\l=C_\l\zeta_\l A_\l,\quad
T_{w_0}A_\l^{-1}=A_\l T_{w_0},\quad
T_{w_0}C_\l\zeta_\l=C_\l\zeta_\l T_{w_0}^{-1}.\leqno(5.8.3)$$
Thus $\psi_{Z_\l}$ is an idempotent.
From Proposition 5.7.1 and \cite{13, \S 37} we get also
$$B_\l C_\l\G_\l T_{w_0}=C_\l\G_\l T_{w_0}B_\l^{-1},\quad
T_{w_0}C_\l\G_\l=C_\l\G_\l T_{w_0}^{-1}.\leqno(5.8.4)$$
Thus $\b_{Z_\l}$ is an idempotent.
Claim 2 is immediate.
\quad\qed
\enddemo

\head 6. The metric and the involution on standard modules\endhead
\subhead 6.1\endsubhead
Set 
$a_{\l\a}=(-q)^{(\rho,\a)},$
$b_{\l\a}=(-q)^{(\rho,\a-2\l)}q^{-d_{\l\a}/2}.$
Let $a_\l$, $b_\l$ be the automorphisms 
of the $\Rb^{\tilde H_\l}$-module $\Wb_{H_\l}$
(resp. of $\Wb'_{H_\l}$) such that
$a_\l(x)=a_{\l\a}x,$ $b_\l(x)=b_{\l\a}x$
for any element $x\in\Wb_{\l\a}$
(resp. $x\in\Wb'_{\l\a}$).
Using (5.8.1) we get
$$b_\l(u\cdot x)=B(u)\cdot b_\l(x),\quad
a_\l(u\cdot x)=A(u)\cdot a_\l(x).$$
Let $\b_\l$, $\b'_\l$ be the automorphisms of 
$\Wb_{H_\l},$ $\Wb'_{H_\l}$ respectively such that
$$\b_\l=T_{w_0}b_\l c_\l\g_\l,\quad
\b'_\l=T_{w_0}b_\l c_\l\g'_\l.$$

\proclaim{Proposition}
\roster
\item We have $\b_\l(u\star x)=\b_{Z_\l}(u)\star\b_\l(x)$
for any $u\in\Ab_\l$, $x\in\Wb_{H_\l}$.
	
\item We have
$\b'_\l(u\star x)=q^{d_{\l\a}-d_{\l\a'}}\b_{Z_\l}(u)\star\b'_\l(x),$
for any 
$u\in\Kb^{\tilde H_\l}(Z_{\l\a\a'})$,
$x\in\Wb'_{H_\l,\a'}$.

\item The maps $\b_\l$, $\b'_\l$ are $\bar{\ }$-semilinear.
Moreover we have $\b_\l^2=\Id$, ${\b'_\l}^2=\Id$.

\item We have $\b_\l(1_\l)=1_\l$, $\b'_\l(1'_\l)=1'_\l$.
\endroster
\endproclaim

\demo{Proof}
Claim 1 follows from Lemma 5.3.3. 
Claim 2 follows from Lemma 5.3.4.
Claim 3 follows from Corollary 5.8.1
and the equality $F_{\l\a}=Z_{\l\a0}$.
Using Lemma 5.7.2 we get
$$\g_\l(1_\l)=1_{\l\nu},\quad
\g'_\l(1'_\l)=1'_{\l\nu},\quad
b_{\l\nu}=1,\quad c_{\l\nu}=r_\l^{-1}1'_{\l\nu}.$$
Thus Claim 4 follows from Proposition 5.6.3.
\quad\qed
\enddemo

\remark{Remark} 
For any closed subgroup $H_\l'\subset H_\l$, the forgetful maps
$\Wb_{H_\l}\to\Wb_{H'_\l}$, $\Wb'_{H_\l}\to\Wb'_{H'_\l}$
commute with the involutions $\b_\l$, $\b'_\l$.
\endremark
\subhead 6.2\endsubhead
For any $\AA$-module $M$, let 
$\hat M$ be the set of formal series in $q^{-1}$
with coefficients in $M$.
We get (see 5.5)
$$\hat\Wb_{H_\l}=\Wb_{H_\l}\otimes_{\Rb^{\tilde H_\l}}\hat\Rb^{\tilde H_\l},\quad
\Wb'_{H_\l}=\Wb'_{H_\l}\otimes_{\Rb^{\tilde H_\l}}\hat\Rb^{\tilde H_\l}.$$
Recall that if $\l=\l_1+\l_2$ in $P^+$ 
then the direct sum of representations of the quiver $\Pi^e$
gives an embedding $\varpi\,:\,Q_{\l_1}\times Q_{\l_2}\hookrightarrow Q_\l$.
Fix a pair of ring isomorphisms 
$$\Rb^{T_{\l_1}}\simeq\ZZ[x_1^{\pm 1},...,x_{\ell_1}^{\pm 1}],\quad
\Rb^{T_{\l_2}}\simeq\ZZ[y_1^{\pm 1},...,y_{\ell_2}^{\pm 1}].$$
We have $\Rb^{T_\l}\simeq\Rb^{T_{\l_1}}\otimes\Rb^{T_{\l_2}}$.
Set
$$\hat\Rb_{\l_1/\l_2}=\ZZ[[q^{-1},y_i/x_j;i,j]]
\otimes_{\ZZ[q^{-1},y_i/x_j;i,j]}\Rb^{\tilde T_\l}.$$
where $1\leq i\leq\ell_1$ and $1\leq j\leq\ell_2$.
Recall that $\kappa$ is the closed embedding $F_\l\hookrightarrow Q_\l$.

\proclaim{Lemma}
\roster
\item The direct image map $\kappa_*$
is an isomorphism $\hat\Wb_{H_\l}\simto\hat\Wb'_{H_\l}$.

\item Assume that $\l=\l_1+\l_2$ in $P^+$. 
Then, there is a unique isomorphism of 
$\hat\Rb_{\l_1/\l_2}\otimes\Ub$-modules
$$\varpi_{\l_1/\l_2}\,:\,
\hat\Rb_{\l_1/\l_2}\otimes_{\Rb^{\tilde T_\l}}
(\Wb_{T_{\l_1}}\otimes_\AA\Wb_{T_{\l_2}})\simto
\hat\Rb_{\l_1/\l_2}\otimes_{\Rb^{\tilde T_\l}}\Wb_{T_\l}$$
such that $1_{\l_1}\otimes 1_{\l_2}\mapsto 1_\l.$
\endroster
\endproclaim

\demo{Proof}
Let first prove Claim 1.
Assume that $H_\l=T_\l$.
We set $\ell=\sum_i\ell_i$.
There is an isomorphism of rings
$\Rb^{T_\l}\simeq\ZZ[z_1^{\pm 1},...,z_\ell^{\pm 1}]$.
Fix $\Rb^{\tilde T_\l}$-bases in $\Wb_\l$, $\Wb'_\l$.
By Thomason's concentration theorem in equivariant $K$-theory and
by \cite{20, Proposition 4.2.2}, the determinant of the 
map $\kappa_*$ in those bases belongs to the set
$$\bigl(\Rb^{\tilde T_\l}\bigr)^\times\cdot\prod_k(1-q^{n_k}z_{i_k}/z_{j_k})$$
for some $i_k,j_k\in[1,\ell]$, $n_k\in\ZZ\setminus\{0\}$.
We can assume that $n_k<0$ for all $k$.
Thus this determinant is invertible in the ring $\hat\Rb^{\tilde T_\l}$.
The case of a general group $H_\l$ follows from Lemma 5.5.
Let us prove Claim 2.
In \cite{24, Proposition 7.10.$(v)$} we define
an embedding of $\Rb^{T_\l}\otimes\Ub$-modules
$$\Delta_W\,:\,\Wb_{T_\l}\to\Wb_{T_{\l_1}}\otimes_\AA\Wb_{T_{\l_2}}.$$
By \cite{24, Theorem 7.12} the map $\Delta_W$ is an isomorphism
whenever $q,x_j,y_i$ are specialized to non-zero complex numbers
such that $y_i/x_j\notin q^{1+\NN}$ for all $i,j$.
Hence it yields an isomorphism of $\hat\Rb_{\l_1/\l_2}\otimes\Ub$-modules
$$\hat\Rb_{\l_1/\l_2}\otimes_{\Rb^{\tilde T_\l}}
(\Wb_{T_{\l_1}}\otimes_\AA\Wb_{T_{\l_2}})\simto
\hat\Rb_{\l_1/\l_2}\otimes_{\Rb^{\tilde T_\l}}\Wb_{T_\l}.$$
The unicity follows from Lemma 5.5.
\quad\qed
\enddemo
\subhead 6.3\endsubhead
Let $a$ be the map from $Q_\l$ to the point.
We consider the pairing of $\Rb^{\tilde H_\l}$-modules
$$(\ :\ )\,:\,\Wb_{H_\l}\times\Wb'_{H_\l}\to\Rb^{\tilde H_\l}$$
given by $(x:y)=a_*(x\otimes y)$, where $\otimes$ is the tor-product
relative to the smooth variety $Q_\l$.
The pairing $(\ :\ )$ is perfect, see \cite{20}.
Note that, $\Wb'_{H_\l}$ beeing a free $\AA$-module,
there is an embedding $\Wb'_{H_\l}\subset\hat\Wb'_{H_\l}$
Let us consider the pairings
$$(\,||\,)\,:\,\Wb_{H_\l}\times\Wb'_{H_\l}\to\Rb^{\tilde H_\l},$$
$$(\,|\,)\,:\,\Wb_{H_\l}\times\Wb_{H_\l}\to\Rb^{\tilde H_\l},\quad
(\,|\,)'\,:\,\Wb'_{H_\l}\times\Wb'_{H_\l}\to\hat\Rb^{\tilde H_\l}$$
such that 
$$(x||y)=(c_\l^{-1} a_\l x:\o^*T_{w_0}^{-1}(y)),\quad
(x|y)=(x||\kappa_*(y)),\quad (x|y)'=(\kappa_*^{-1}(x)||y),$$
see Lemma 6.2.1.
Let $\partial\,:\,\Rb^{\tilde H_\l}\to\AA$ 
be the group homomorphism such that $\partial(q)=q$,
and $\partial(V)=0$ if
$V$ is a non trivial simple $H_\l$-module.

\proclaim{Proposition}
\roster 
\item We have $(x|y)=(y|x)^\dagger$.

\item We have $(x 1_\l|y 1_\l)=xy^\dagger$, 
for all $x,y\in\Rb^{\tilde H_\l}$.

\item We have $(u\cdot x|y)=(x|\psi(u)\cdot y).$

\item The pairing $(\,|\,)$ is uniquely determined by conditions {\rm 2}
and {\rm 3}.

\item We have $(\b_\l(x)||y)=\overline{(x||\b'_\l(y))}$.

\item The pairing of $\AA$-modules $\partial(\,||\,)$ is perfect.

\item Claims {\rm 1} and {\rm 2} hold for the pairing $(\,|\,)'$ also.
\endroster
\endproclaim

\demo{Proof}
First, note that $a_{\l0}=1$.
By Lemma 5.7.2 we have $c_{\l0}=(r_\l^\dagger)^{-1} 1'_\l$.
By Proposition 5.6.3 we have 
$\o^*T_{w_0}^{-1}(1_\l)=(r_\l^\dagger)^{-1} 1_\l$.
Thus
$$\matrix
(x1_\l|y1_\l)
&=(c_{\l0}^{-1} a_{\l0}x1_\l:\o^*T_{w_0}^{-1}(y1_\l)),\hfill\cr
&=(x1_\l:y^\dagger 1_\l),\hfill\cr
&=xy^\dagger.\hfill\cr
\endmatrix$$
Claim 2 is proved.
Fix $u\in\Ub,$ 
$x\in\Wb_{H_\l,\a}$, 
$y\in\Wb'_{H_\l,\a}$.
For all $w\in W$ let $1_{\l,w*0}\boxtimes x\in\Kb^{\tilde H_\l}(Z_{\l,w*0,\a})$
be the obvious element. 
We have
$$(1_{\l,w*0}\boxtimes x)\star y=(x:y)1_{\l,w*0}.$$
Thus, the associativity of $\star$ gives
$$(x\cdot u:y)1_\l=(1_\l\boxtimes x)\star\Phi_\l(u)\star y=(x:u\cdot y)1_\l.$$
Thus we get
$$(u\cdot x:y)=(x:\phi^*\Phi_\l(u)\star y).\leqno(6.3.0)$$
Assume now that
$x\in\Wb_{H_\l,\a'}$, 
$y\in\Wb_{H_\l,\a}$.
Using (6.3.0), Lemma 5.2.3, Proposition 5.6.1,
(5.7.6) and (5.8.3) we get
$$\matrix
(u\cdot x|y)
&=
(c_{\l\a}^{-1}a_{\l\a}\Phi_\l(u)\star x:\o^* T_{w_0}^{-1}\kappa_*(y))
\hfill\cr
&=
(c_{\l\a'}^{-1}a_{\l\a'} x:\phi^*C_\l^{-1} A_\l\Phi_\l(u)\star
\o^* T_{w_0}^{-1}\kappa_*(y))
\hfill\cr
&=
(x|T_{w_0}\zeta_\l C_\l^{-1}A_\l\Phi_\l(u)\star y)\hfill\cr
&=
(x|T_{w_0}A_\l C_\l\zeta_\l\Phi_\l(u)\star y)\hfill\cr
&=
(x|\psi_{Z_\l}\Phi_\l(u)\star y).\hfill
\endmatrix$$
Then, apply Corollary 5.8.2.
Claim 3 is proved.
Claim 4 follows from Lemma 5.5.2.
Assume now that
$x\in\Wb_{H_\l,\a}$, 
$y\in\Wb'_{H_\l,\a}$ as above.
For all $w$ we have
$$1_{\l,w*0}\star(1_{\l,w*0}\boxtimes x)=x.$$
Fix $u_\l$ such that
$\check T_{w_0}(1_{\l\nu})=u_\l 1_\l.$
Using Proposition 5.6 we get
$$u_\l 1_\l\star T_{w_0}^{-1}(1_{\l\nu}\boxtimes x)=\check T_{w_0}(x),\quad
T_{w_0}^{-1}(1_{\l\nu}\boxtimes x)\star T^{-1}_{w_0}(y)=s_\l^{-1}(x:y)1_\l,$$
i.e.
$$\matrix
T_{w_0}^{-1}(1_{\l\nu}\boxtimes x)=u_\l^{-1}1_\l\boxtimes\check T_{w_0}(x),
\hfill\cr\cr
(1_\l\boxtimes\check T_{w_0}(x))\star T^{-1}_{w_0}(y)=
u_\l s_\l^{-1}(x:y)1_\l.\hfill
\endmatrix
\leqno(6.3.1)$$
This yields
$$(\check T_{w_0}(x):T^{-1}_{w_0}(y))=u_\l s_\l^{-1}(x:y).\leqno(6.3.2)$$
Claim 5 is analoguous to \cite{14, Lemma 12.15}.
First, using (4.5.1) one gets
$$\overline{(x:\DD_{Q_{\l\a}}(y))}^\dagger=
q^{d_{\l\a}}(\DD_{F_{\l\a}}(x):y).\leqno(6.3.3)$$
Note that
$$b_{\l,w_0*\a}a_{\l\a}=q^{-d_{\l\a}/2}.\leqno(6.3.4)$$
Using (6.3.0), Lemma 5.7.2, (6.3.4) and (6.3.3) we get
$$\matrix
\overline{(x||\b'_\l(y))}&=
q^{-3d_{\l\a}/2}
\overline{(c_{\l\a}^{-1}a_{\l\a} x:
c_{\l\a} b_{\l,w_0*\a}\DD_{Q_{\l\a}}(y))}\hfill\cr
&=q^{-d_{\l\a}}\overline{(x:\DD_{Q_{\l\a}}(y))}\hfill\cr
&=(\DD_{F_{\l\a}}(x):y)^\dagger.\hfill
\endmatrix$$
Using Lemma 5.7.2 and (6.3.1) we get 
$$T_{w_0}^{-1}C_\l\zeta_\l(x\boxtimes 1_\l)=
T_{w_0}^{-1}C_\l(1_{\l\nu}\boxtimes\o^*x)=
u_\l^{-1}r_\l^{-1}1_\l\boxtimes\check T_{w_0}c_\l^{-1}\o^*x.$$
The same argument as for (6.3.1) gives
$$T_{w_0}(x\boxtimes 1_\l)=s_\l^{-1}T_{w_0}(x)\boxtimes 1_{\l\nu}.$$
Recall that $\o$ is $\dagger$-linear, $(r_\l s_\l)^\dagger=r_\l s_\l$, and
$c_{\l0}=(r_\l^\dagger)^{-1}1'_\l$.
Thus we get
$$C_\l\zeta_\l T_{w_0}(x\boxtimes 1_\l)=
r_\l^{-1}s_\l^{-1}(1_\l\boxtimes c_\l^{-1}\o^*T_{w_0}x).$$
Thus, (5.8.3) gives
$\check T_{w_0} c_\l^{-1}\o^*(x)=u_\l s_\l^{-1}c_\l^{-1}\o^*T_{w_0}(x),$
i.e. 
$$T_{w_0} c_\l\o^*(x)=c_\l\o^*\check T_{w_0}(u_\l^{-1}s_\l x)\leqno(6.3.5)$$
for all $x\in\Wb_{H_\l}$. 
Using (6.3.4), (6.3.5), (6.3.2) we get
$$\matrix
(\b_\l(x)||y)
&=q^{d_{\l\a}/2}
(c_{\l\a}^{-1}a_{\l\a}T_{w_0}b_{\l,w_0*\a}c_{\l,w_0*\a}\o^*\DD_{F_{\l\a}}(x):
\o^*T_{w_0}^{-1}(y))
\hfill\cr
&=(T_{w_0}c_{\l,w_0*\a}\o^*\DD_{F_{\l\a}}(x):c_{\l\a}^{-1}\o^*T_{w_0}^{-1}(y))
\hfill\cr
&=(u_\l^{-1}s_\l\check T_{w_0}\DD_{F_{\l\a}}(x):T_{w_0}^{-1}(y))^\dagger
\hfill\cr
&=
(\DD_{F_{\l\a}}(x):y)^\dagger.\hfill
\endmatrix$$
Claim 5 is proved.
Claim 1 follows from (6.3.2), (6.3.5) and Lemma 5.7.2.
Indeed
$$\matrix
(x|y)
&=(c_\l^{-1}a_\l x:\o^*T_{w_0}^{-1}\kappa_*(y))
\hfill\cr
&=(u_\l s_\l^{-1}\check T_{w_0}^{-1}\o^*c_\l^{-1}a_\l\kappa_*(x):y)^\dagger
\hfill\cr
&=(\o^*c_\l^{-1}T_{w_0}^{-1}a_\l\kappa_*(x):y)^\dagger
\hfill\cr
&=(y|x)^\dagger.\hfill
\endmatrix$$
Claim 6 follows from the Schur Lemma and the fact that
$(\, :\,)$ is a perfect pairing of $\Rb^{\tilde H_\l}$-modules.
\quad\qed
\enddemo

\remark{Remarks}
\roster
\item The pairings $(\,|\,)$, $(\,|\,)'$ 
are obviously compatible with the forgetful
maps, see Remark 6.1.

\item The $\Ub$-module $\Wb'_{H_\l}$ has the following algebraic
interpretation : 
let $\Wb_{H_\l}^*$ be $\Wb_{H_\l}'$ with the new action of $\Ub$, 
denoted by $\diamond$, such that
$u\diamond x=\phi^*\Phi_\l S(u)\star x,$
where $S$ is the antipode;
then, the $\Ub$-module $\Wb^*_{H_\l}$ is the right dual of $\Wb_{H_\l}$.
\endroster
\endremark

\head 7. Construction of the signed basis\endhead
\subhead 7.1\endsubhead
Following Lusztig we consider the sets
$$\Bc'_{H_\l}=\{\bb\in\Wb'_{H_\l}\,|\,\b'_\l(\bb)=\bb,\,
\partial(\bb|\bb)'\in 1+q^{-1}\ZZ[[q^{-1}]]\}, $$
$$\Bc_{H_\l}=\{\bb\in\Wb_{H_\l}\,|\,\b_\l(\bb)=\bb,\,
\partial(\bb|\bb)\in 1+q^{-1}\ZZ[q^{-1}]\}. $$
We also set $\Bc'_\l=\Bc'_{G_\l}$, $\Bc_\l=\Bc_{G_\l}$.

\proclaim{Proposition}
\roster
\item If the subset $\Bb\subset\Bc_{H_\l}$ satisfies:

\itemitem - $\Bb$ is a basis of the $\AA$-module $\Wb_{H_\l}$,

\itemitem - for any elements $\bb,\bb'\in\Bb$ we have
$\partial(\bb|\bb')\in\delta_{\bb,\bb'}+q^{-1}\ZZ[q^{-1}]$.

\noindent
Then $\Bc_{H_\l}=\pm\Bb$.
The similar statement holds for $\Bc'_{H_\l}$.

\item We have $x1_\l,xr_\l 1_{\l\nu}\in\Bc_{H_\l}$,
and $x1'_\l,xr_\l 1'_{\l\nu}\in\Bc'_{H_\l}$,
for any $x\in\Xb^{H_\l}$.
\endroster
\endproclaim

\demo{Proof}
Claim 1 is standard, see \cite{14, \S 12.20} for instance. 
We reproduce a proof here for the convenience of the reader.
Fix an element $\bb\in\Bc_{H_\l}$. Set $\bb=\sum_ip_i\bb_i$ where
$\bb_i\in\Bb$ and $p_i\in\AA$. Fix $n\in\ZZ$ such that 
$p_i\in q^n\ZZ[q^{-1}]$ for all $i$ and 
$p_i\notin q^{n-1}\ZZ[q^{-1}]$ for some $i$.
For all $i$ let $p_{in}\in\ZZ$ be such that $p_i\in p_{in}q^n+q^{n-1}\ZZ[q^{-1}]$.
Then, $\sum_ip_{in}^2>0$. Thus,
$$\partial(\bb|\bb)\in q^{2n}\Sum_ip_{in}^2+q^{2n-1}\ZZ[q^{-1}].$$
On the other hand, we have $\partial(\bb|\bb)\in 1+q^{-1}\ZZ[q^{-1}]$.
It follows that $n=0$ and $\sum_ip_{in}^2=1$.
Since $\b_\l(\bb)=\bb$ and $\b_\l(\bb_i)=\bb_i$ for all $i$,
we must have $\bar p_i=p_i$ for all $i$.
Hence $p_i\in\ZZ$ for all $i$.
Then $\sum_ip_i^2=1$. Thus $\bb\in\pm\Bb$.
Let us prove Claim 2.
By Proposition 6.1.4 and 6.3.2 we have
$x1_\l\in\Bc_{H_\l}$, $x1'_\l\in\Bc'_{H_\l}$.
Hence, using \cite{11} we get 
$T_{w_0}^{-1}(x1_\l)\in\Bc_{H_\l}$, 
$T_{w_0}^{-1}(x1'_\l)\in\Bc'_{H_\l}$.
Finally, Proposition 5.6.3 gives
$T_{w_0}^{-1}(x1_\l)=xr_\l1_{\l\nu}.$
We are done.
\quad\qed
\enddemo

\remark{Remark}
In general $1_{\l\a}\notin\Bc_{H_\l}$.
\endremark
\subhead 7.2\endsubhead
For any $\l\in P^+$ let $V(\l)$ be Kashiwara's maximal integrable module.
By definition, $V(\l)$ is the free $\AA$-module with the action
of the algebra $\Ub$ such that there is a weight vector
$v_\l$ of weight $\l$ which generates $V(\l)$ 
and satisfies the following defining relations :
$$\matrix
\Ub_\a(v_\l)=0\ \text{for any}\ 
\a\in Q\setminus\{0\}\ \text{s.t.}\ 
(\a,\l)\geq 0,\cr\cr
\fb_i^{1+\ell_i}(v_\l)=0\ \text{if}\ i\neq 0,\quad
\eb_0^{1+(\theta,\l)}(v_\l)=0,
\endmatrix\leqno(7.2.1)$$
see \cite{9, \S 5.1}.
It is proved in \cite{8} that the module $V(\l)$ admits a global basis.
Let $\Bb(\l)$ be this basis. 
The element $v_\l$ belongs to $\Bb(\l)$.
Let $\bar{}\,:\,V(\l)\to V(\l)$ be the unique
$\AA$-antilinear map such that $\bar\bb=\bb$ for all elements
$\bb\in\Bb(\l)$.
It is conjectured in \cite{24, Remark 7.19} that 
there is an isomorphism of $\Ub$-modules $V(\l)\to\Wb_\l$ 
such that $v_\l\mapsto 1_\l$.
Let us first consider the case $\l=\o_i$.
Let $W(\o_i)'$ be the fundamental simple finite dimensional
$\Ub'$-module associated to the weight $\o_i$,
see \cite{9, (5.7)}, \cite{1,\S 1.3}.
Let $W(\o_i)\subset W(\o_i)'$ be the corresponding $\AA$-form
For any $\Ub$-module $M$ and any formal variable $z$, let $M_z$ be the
representation of $\Ub$ on the space $M[z^{\pm 1}]$ such that
$(\xb_{jr}^\pm)^{(n)}\mapsto(\xb_{jr}^\pm)^{(n)}\otimes z^{r},$
$\kb^\pm_{jr}\mapsto\kb^\pm_{jr}\otimes z^r.$
Fix a weight vector $w_{\o_i}\in W(\o_i)$ of weight $\o_i$.
By \cite{9, Theorem 5.15.$(viii)$} 
there is a unique  isomorphism of $\Ub$-modules 
$$V(\o_i)\simto W(\o_i)_z\leqno(7.2.2)$$ 
such that $v_{\o_i}\mapsto w_{\o_i}.$
The product by $z$ is an automorphism of $\Ub$-modules.
It preserves the basis $\Bb(\o_i)$.
There is a unique basis $\Bb^0(\o_i)$ of $W(\o_i)$ 
such that the map (7.2.2) takes 
$\Bb(\o_i)$ to $\bigsqcup_{n\in\ZZ}z^n\Bb^0(\o_i)$,
see \cite{9, Theorem 5.15.$(iii)$}.

The group $G_{\o_i}$ being isomorphic to $\CC^\times$ we 
identify $\Rb^{G_{\o_i}}$ with $\ZZ[z_i^{\pm 1}]$ in the usual way. 

\proclaim{Theorem A}
\roster 
\item There is a unique element $a_i\in\QQ(q)^\times$
and a unique isomorphism of $\Ub$-modules 
$\phi\,:\,V(\o_i)\simto\Wb_{\o_i}$ such that $v_{\o_i}\mapsto 1_{\o_i}$ and
the multiplication by $z$ is mapped to the multiplication by $a_iz_i.$
Moreover $a_i=\pm 1$.
\footnote
{H. Nakajima remarked that
Theorem A.1 was not stated correctly in a previous version of the paper.
He also mentioned to us that the isomorphism of $\Ub$-modules
$V(\o_i)\simto\Wb_{\o_i}$ was known to him.}

\item 
Assume that
$\la\,|\,\ra\,:\,V(\o_i)\times V(\o_i)\to\AA$ is a symmetric perfect pairing
of $\AA$-modules such that
$\la z^nv_{\o_i}|z^mv_{\o_i}\ra=\delta_{n,m}$
and
$\la u\cdot x|y\ra=\la x|\psi(u)\cdot y\ra.$
Then
$$\pm\Bb(\o_i)=\{\bb\in V(\o_i)\,|\,\bar\bb=\bb,\,
\la\bb|\bb\ra\in 1+q^{-1}\ZZ[q^{-1}]\}.$$
Moreover, $\la\bb|\bb'\ra\in q^{-1}\ZZ[q^{-1}]$ if $\bb,\bb'\in\Bb(\o_i)$
and $\bb\neq\bb'$.

\item $\Bc_{G_{\o_i}}=\pm\phi(\Bb(\o_i)).$
It is a signed basis of $\Wb_{\o_i}$.
\endroster
\endproclaim

\demo{Proof}
Let us prove Claim 1.
We identify $\Kb^{\CC^\times}(F_{\o_i})$ with the specialization
of the $\Ub$-module $\Wb_{\o_i}$ at the maximal ideal of $\Rb^{G_{\o_i}}$
associated to $1\in G_{\o_i}$.
There is a unique element $a_i\in\QQ(q)^\times$ and a unique
isomorphism of $\Ub'$-modules
$\QQ(q)\otimes_\AA\Kb^{\CC^\times}(F_{\o_i})\simto W(\o_i)'_{a_i}$ 
such that $1_{\o_i}\mapsto w_{\o_i}$,
since both $\Ub'$-modules are simple, see \cite{20}.
The $\Ub$-modules $\Kb^{\CC^\times}(F_{\o_i})$, $W(\o_i)_{a_i}$ 
beeing cyclic generated by $1_{\o_i}$, $w_{\o_i}$,
we get an isomorphism
$\Kb^{\CC^\times}(F_{\o_i})\simto W(\o_i)_{a_i}$ 
such that $1_{\o_i}\mapsto w_{\o_i}$.
The identification of the group $G_{\o_i}$  with $\CC^\times$ is such that
for any $(B,p,q)\in M_{\o_i\a}$ and any $g_{\o_i}\in G_{\o_i}$
we have
$$(1,g_{\o_i},1)\cdot(B,p,q)=(B,g^{-1}_{\o_i}p,g_{\o_i}q).$$
Since $p_j=0$, $q_j=0$ if $j\neq i$ we have 
$$(1,g_{\o_i},1)\cdot(B,p,q)=(1,1,g_\a)\cdot(B,p,q)$$
for $g_\a=(g_{a_j})_j$ with $g_{a_j}=g_{\o_i}^{-1}\Id_{\CC^{a_j}}.$
Then the group $G_{\o_i}$ acts trivially on $F_{\o_i}$ and the 
natural isomorphism of $\AA[z_i^{\pm 1}]$-modules
$$\Wb_{\o_i}=\Kb^{G_{\o_i}\times\CC^\times}(F_{\o_i})\simto 
\Kb^{\CC^\times}(F_{\o_i})[z_i^{\pm 1}]$$
takes $\Vc_j$ to $\Vc_j\otimes z_i$, and $\Wc_j$ to $\Wc_j\otimes z_i$.
In particular we have
$$
\Wedge_{\Vc_\a}\mapsto\Wedge_{\Vc_\a}\otimes z_i^{(\rho,\a)},\ 
\Wedge_{\Fc_{j;\a}^+}\mapsto\Wedge_{\Fc_{j;\a}^+}\otimes z_i^{f_{j;\a}^+},\ 
\Wedge_{\Fc_{j;\a}^-}\mapsto\Wedge_{\Fc_{j;\a}^-}\otimes z_i^{f_{j;\a}^-}
$$
and
$\xb_{jr}^\pm\mapsto\xb_{jr}^\pm\otimes z_i^{r},$
$\kb_{jr}^\pm\mapsto\kb_{jr}^\pm\otimes z_i^r,$
because $t_{\a+\a_j}-t_\a=f^+_{j;\a}-f^-_{j;\a}.$ 
Hence $\Wb_{\o_i}\simeq(W(\o_i)_{a_i})_{z_i}\simeq W(\o_i)_{a_iz_i}$.

The map $\phi$ takes the involution $v\mapsto\bar v$ on $V(\o_i)$ to the 
involution $\beta_{\o_i}$ on $\Wb_{\o_i}$ since both $\Ub$-modules are cyclic
and the involutions are compatible with $u\mapsto\bar u$ on $\Ub$.
The map $\beta_{\o_i}$ is $z_i$ linear by Corollary 5.8.1,
the map $v\mapsto\bar v$ on $V(\o_i)$ is $z$ linear since product by $z$
preserves $\Bb(\o_i)$. Hence $\bar a_i=a_i$.
The $j$-th Drinfeld polynomial
of $\Kb^{\CC^\times}(F_{\o_i})$ is $P_j(t)=(t-q^{-\cb})^{\delta_{ij}}$. 
Hence the $j$-th Drinfeld polynomial
of $W(\o_i)$ is $P_j(t)=(t-a_i^{-1}q^{-\cb})^{\delta_{ij}}$. 
An easy computation shows that the elements $\hb_{i,\pm 1}$ act
on the vector $w_{\o_i}\in W(\o_i)$ as follows
$\hb_{i,\pm 1}(w_{\o_i})=a_i^{\pm 1}q^{\pm\cb}w_{\o_i}$.
Since $\hb_{i,\pm 1}$ belongs to $\Ub$ the element 
$a_iq^{-\cb}$ belongs to $\AA$ and is invertible.
Thus $a_i=\pm 1$ because $a_i\in\pm q^\ZZ$ and $a_i=\bar a_i$.

Let us prove Claim 2.
We first recall some well-known fact.
Let $\dot\Bb$ be the canonical basis of $\dot\Ub$, see \cite{13, \S 25.2}.
By \cite{8, \S 8} there is a subset $\Ib(\l)\subset\dot\Bb$ 
such that the space 
$I(\l)=\bigoplus_{\bb\in\Ib(\l)}\AA\,\bb\subset\dot\Ub\eta_\l$ 
is a left $\dot\Ub$-submodule,
and such that there is a unique isomorphism of $\dot\Ub$-modules
$\dot\Ub\eta_\l/I(\l)\to V(\l)$ which takes $\eta_\l$ to $v_\l$.
The $\Ub$-module $V(\l)$ beeing integrable, 
Kashiwara's modified operators
$\tilde e_j$, $\tilde f_j$, $j\in I\cup\{0\}$, act on $V(\l)$.
Let $L(\l)\subset V(\l)$ be the $\ZZ[q^{-1}]$-lattice 
linearly spanned by $\Bb(\l)$.
It is stable by the operators $\tilde e_j$, $\tilde f_j$, 
see \cite{8, Proposition 9.1}, and contains the element $v_{\l}$.
The induced operators on the quotient $L(\l)/q^{-1}L(\l)$ 
are still denoted by $\tilde e_j$, $\tilde f_j$. 
Let $\pi\,:\,L(\l)\to L(\l)/q^{-1}L(\l)$ be the projection.
We set $B(\l)=\pi\bigl(\Bb(\l)\bigr)$.
It is known that $\tilde e_j$, $\tilde f_j$ 
take $B(\l)$ to $B(\l)\sqcup\{0\}$.

We now assume that $\l=\o_i$.
Then the operators $\tilde e_j$, $\tilde f_j$ are $z$-linear. 
Let $L^0(\o_i)$ be the $\ZZ[q^{-1}]$-module spanned by $\Bb^0(\o_i)$.
Let $B^0(\o_i)$ be the projection of $\Bb^0(\o_i)$ in 
$L^0(\o_i)/q^{-1}L^0(\o_i)$.
There is an isomorphism of crystals 
$$\bigl(B^0(\o_i),L^0(\o_i)\bigr)\simeq\bigl(B(\o_i),L(\o_i)\bigr)/(z-1).$$
Any element in $B^0(\o_i)$ can be reached at $w_{\o_i}$
after applying a monomial in the operators
$\tilde e_j$, $j\in I\cup\{0\}$,
see \cite{1, Lemma 1.5.(1) and (2)} 
and \cite{9, Proposition 5.4.$(i)$}.
Thus any element in $B(\o_i)$ can be reached at $\{z^mv_{\o_i}\, ;\,m\in\ZZ\}$
after applying a monomial in the operators
$\tilde e_j$, $j\in I\cup\{0\}$.

Set $L(\o_i)^\infty=\bigcup_{k\geq 0}L(\o_i)^{k},$ where
$$L(\o_i)^{k}=\sum_{\ell\leq k}
\sum_{j_1,...j_\ell}\ZZ[q^{-1}, z^{\pm 1}]\,
\tilde f_{j_1}\cdots\tilde f_{j_\ell}(v_{\o_i}).$$
We claim that
$$L(\o_i)=L(\o_i)^\infty+q^{-1} L(\o_i),\leqno(7.2.3)$$
$$\la L(\o_i)\,|\,L(\o_i)\ra\subseteq\ZZ[q^{-1}],\quad
\la\tilde f_j(x)|y\ra\in\la x|\tilde e_j(y)\ra+q^{-1}\ZZ[q^{-1}],
\leqno(7.2.4)$$
$\forall j\in I\cup\{0\}$, $\forall x,y\in L(\o_i).$
Claim (7.2.3) is obvious.
To prove (7.2.4) we use the following Lemma, 
whose proof is given after the proof of the proposition.

\proclaim{Lemma}
Fix $j\in I\cup\{0\}$. 
For any $x\in V(\o_i)$ fix elements $x_r\in V(\o_i)$, $r\in [0,t]$, such that 
$x=\sum_{r=0}^t\fb_j^{(r)}(x_r)$ and $\eb_j(x_r)=0$.
\roster
\item If $x\in L(\o_i)$ then $x_r\in L(\o_i)$.

\item If $x\in\Bb(\o_i)$ then there is $r_0\in[0,t]$ such that
$x_{r_0}\in\Bb(\o_i)+q^{-1}L(\o_i)$ and
$x_r\in q^{-1}L(\o_i)$ if $r\neq r_0$.

\item Fix $\l,\mu\in P$.
Assume that $x\in L(\o_i)_\l$,
and $\la x\,|\,L(\o_i)\ra\subseteq\ZZ[q^{-1}]$.
Fix $y\in L(\o_i)_\mu$, with $\mu=\l-\a_j$ and fix elements $y_s\in V(\o_i)$, 
$s\in [0,u]$, 
such that $y=\sum_{s=0}^u\fb_j^{(s)}(y_s)$ and $\eb_j(y_s)=0$.
Then $\la x_r\,|\,y_s\ra\in\ZZ[q^{-1}]$ for all $r,s$.
\endroster
\endproclaim

The $\Ub$-module $V(\o_i)$ is endowed with its $\hat P$-gradation, see \cite{9}.
For any $\hat\mu\in\hat P$ let $V(\o_i)_{\hat\mu}\subset V(\o_i)$,
$L(\o_i)_{\hat\mu}=L(\o_i)\cap V(\o_i)_{\hat\mu}$
and $L(\o_i)_{\hat\mu}^{k}=L(\o_i)^k\cap V(\o_i)_{\hat\mu}$
be the corresponding weight subspaces. 
Then
$$\matrix
\la V(\o_i)_{\hat\mu_1}\,|\,V(\o_i)_{\hat\mu_2}\ra\neq 0
\,\Rightarrow\,\hat\mu_1-\hat\mu_2\in\ZZ\delta,\hfill\cr\cr
L(\o_i)_{\hat\mu}^{k}=
\Sum_j\tilde f_j\bigl(L(\o_i)_{\hat\mu+\a_j}^{k-1}\bigr),\hfill\cr\cr
\tilde e_j\bigl(L(\o_i)_{\hat\mu}\bigr)\subseteq
L(\o_i)_{\hat\mu+\a_j}.\hfill
\endmatrix\leqno(7.2.5)$$

We first prove by induction on $k$ that
$$\la L(\o_i)^{k}\,|\,L(\o_i)\ra\subseteq\ZZ[q^{-1}],\leqno(7.2.6)$$
$$\la\tilde f_j(x)|y\ra\in\la x|\tilde e_j(y)\ra+q^{-1}\ZZ[q^{-1}],
\leqno(7.2.7)$$
$\forall j\in I\cup\{0\}$, $\forall x\in L(\o_i)^{k},$ 
$\forall y\in L(\o_i).$
We have $L(\o_i)^0=\ZZ[q^{-1},z^{\pm 1}]v_{\o_i}$,
and $L(\o_i)_{\o_i+n\delta}=\ZZ[q^{-1}] z^nv_{\o_i}$.
Thus, (7.2.6) for $k=0$ reduces to
$$\la\ZZ[q^{-1}]z^mv_{\o_i}\,|\,\ZZ[q^{-1}]z^nv_{\o_i}\ra\subseteq\ZZ[q^{-1}],$$
which follows from
$\la z^mv_{\o_i}\,|\,z^nv_{\o_i}\ra=\delta_{mn}$.
Similarly, (7.2.7) for $k=0$ reduces to
$$\la\tilde f_j(z^nv_{\o_i})|y\ra\in
\la z^nv_{\o_i}|\tilde e_j(y)\ra+q^{-1}\ZZ[q^{-1}].\leqno(7.2.8)$$
This is obvious if $j=0$  because
$\la\tilde f_j(z^nv_{\o_i})|y\ra\neq 0$ or
$\la z^nv_{\o_i}|\tilde e_j(y)\ra\neq 0$ implies that
$y\in\bigoplus_nV(\o_i)_{\o_i+\theta+n\delta}$,
and $V(\o_i)_{\o_i+\theta+n\delta}=\{0\}$ for all $n$,
see \cite{9, Proposition 5.14.$(i)$} for instance.
If $j\neq 0$, (7.2.8) is proved as follows.
Since $\eb_j(z^nv_{\o_i})=0$, we have
$\tilde f_j(z^nv_{\o_i})=\fb_j(z^nv_{\o_i})$.
Fix elements $y_s\in V(\o_i)$, $s\in[0,u]$, as in Lemma 7.2.3.
Then, $y_s\in L(\o_i)$ by Lemma 7.2.1.
We must show that
$$\la\fb_j(z^nv_{\o_i})\,|\,\Sum_{s\geq 0}\fb_j^{(s)}(y_s)\ra\in
\la z^nv_{\o_i}\,|\,\Sum_{s\geq 1}\fb_j^{(s-1)}(y_s)\ra+q^{-1}\ZZ[q^{-1}].$$
The computation in \cite{13, Proposition 19.1.3} gives the result, since
$\la z^nv_{\o_i}\,|\,y_s\ra\in\ZZ[q^{-1}]$ for all $s$, by (7.2.6) for $k=0$.
We may therefore assume that (7.2.6), (7.2.7)
are already known for $k-1$ with $k>0$.
Using (7.2.5) we see that (7.2.6) for $k$ follows from (7.2.6), (7.2.7) for $k-1$.
Finally, (7.2.7) for $k$ is proved as in \cite{13, Proposition 19.1.3} using
(7.2.6) for $k$, and Lemma 7.2.3.

From (7.2.6) we get
$$\la L(\o_i)^\infty\,|\,L(\o_i)\ra\subseteq\ZZ[q^{-1}].$$
On the other hand there is a positive integer $a$ such that
$$\la L(\o_i)\,|\,L(\o_i)\ra\subseteq q^a\ZZ[q^{-1}]$$
because $\Bb^0(\o_i)$ is finite.
Using (7.2.3) it yields
$$\la L(\o_i)\,|\,L(\o_i)\ra\subseteq\ZZ[q^{-1}].$$
Similarly, given $x,y\in L(\o_i)$ 
we fix $x^\infty\in L(\o_i)^\infty$ such that 
$x-x^\infty\in q^{-1}L(\o_i)$, see (7.2.3).
Then (7.2.7) yields
$$\la\tilde f_j(x^\infty)|y\ra\in\la x^\infty|\tilde e_j(y)\ra+
q^{-1}\ZZ[q^{-1}].$$
Thus
$$\matrix
\la\tilde f_j(x)|y\ra\in
\la\tilde f_j(x^\infty)|y\ra+q^{-1}\ZZ[q^{-1}]&=
\la x^\infty|\tilde e_j(y)\ra+q^{-1}\ZZ[q^{-1}]\cr\cr
&=
\la x|\tilde e_j(y)\ra+q^{-1}\ZZ[q^{-1}].
\endmatrix$$
We have proved (7.2.4).

Then, Claim 2 is proved as in \cite{13, Lemma 19.1.4}, 
using (7.2.3), (7.2.4) and Proposition 7.1.1. 
More precisely for any element $\bb\in\Bb(\o_i)$,
let $\ell(\bb)$ be the smallest $k\geq 0$ such that 
$\bb\in L(\o_i)^k+q^{-1}L(\o_i)$, see (7.2.3).
For any $\bb,\bb'\in\Bb(\o_i)$ 
we prove by induction on $\ell(\bb)$ that
$$\la\bb\,|\,\bb'\ra\in q^{-1}\ZZ[q^{-1}]\ \roman{if}\ \bb\neq\bb',
\leqno(7.2.9)$$
$$\la\bb\,|\,\bb\ra\in 1+q^{-1}\ZZ[q^{-1}].\leqno(7.2.10)$$
If $\ell(\bb)=0$ then $\bb=z^nv_{\o_i}$ for some $n\in\ZZ$.
Thus, if $\la\bb\,|\,\bb'\ra\neq 0$ then $\bb'=z^mv_{\o_i}$ for some $m$
and both statements are obvious.
Fix $k>0$.
Assume that (7.2.9), (7.2.10) hold for any $\bb,\bb'$ such that $\ell(\bb)<k$.
Fix $\bb,\bb'$ such that $\ell(\bb)=k$.
By (7.2.3) there is an integer $j\in I\cup\{0\}$ and
an element $\bb_1\in\Bb(\o_i)$ such that
$\tilde f_j(\bb_1)\in\bb+q^{-1}L(\o_i)$
and $\ell(\bb_1)=k-1$.
Using (7.2.4) we get
$$\la\bb\,|\,\bb'\ra\in
\la\tilde f_j(\bb_1)\,|\,\bb'\ra+q^{-1}\ZZ[q^{-1}]
=\la\bb_1\,|\,\tilde e_j(\bb')\ra+q^{-1}\ZZ[q^{-1}].
$$
We have $\tilde e_j(\bb)\in\bb_1+q^{-1}L(\o_i)$.
Thus,
$$\la\bb\,|\,\bb\ra\in
\la\bb_1\,|\,\bb_1\ra+q^{-1}\ZZ[q^{-1}].$$
Hence (7.2.10) for $k$ follows from (7.2.10) for $k-1$.
If $\bb\neq\bb'$, either $\tilde e_j(\bb')\in q^{-1}L(\o_i)$, 
and then
$$\la\bb\,|\,\bb'\ra\in q^{-1}\ZZ[q^{-1}]$$
by (7.2.4), or there is an element $\bb'_1\in\Bb(\o_i)$
such that $\tilde e_j(\bb')\in\bb'_1+q^{-1}L(\o_i)$.
In the last case $\bb_1\neq\bb'_1$
(else applying $\tilde f_j$ we would get $\bb\in\bb'+q^{-1}L(\o_i)$, 
and thus $\bb=\bb'$).
Hence (7.2.9) for $k$ follows from (7.2.9) for $k-1$.
Finally, Claim 2 follows from 
(7.2.9), (7.2.10) and Proposition 7.1.1. 

Claim 3 is obvious from Claim 1 and Claim 2.
\quad\qed
\enddemo

\demo{Proof of Lemma 7.2}
Claims 1,2 generalize \cite{13, Lemma 18.2.2} 
to the non-highest weight module case. 
The proof follows \cite{13, Lemma 18.2.2}.
Note that we only use Claim 1; 
Claim 2 is given for the sake of completness.

We prove Claim 1 by induction on $t$.
It is obvious if $t=0$. Since
$$\tilde e_j(x)=\sum_{r=0}^{t-1}\fb_j^{(r)}(x_{r+1}),\quad
\tilde e_j(L(\o_i))\subseteq L(\o_i),$$
we get $\sum_{r=0}^{t-1}\fb_j^{(r)}(x_{r+1})\in L(\o_i)$.
The induction hypothesis for $t-1$ gives $x_{r+1}\in L(\o_i)$ for all 
$r\in [0,t-1]$.
Since $\eb_j(x_{r+1})=0$ we have 
$\fb_j^{(r+1)}(x_{r+1})=\tilde f_j^{r+1}(x_{r+1})$.
Since $\tilde f_j(L(\o_i))\subset L(\o_i)$ we get
$\fb_j^{(r+1)}(x_{r+1})\in L(\o_i)$.
Using $x\in L(\o_i)$, $\fb_j^{(r)}(x_r)\in L(\o_i)$
for all $r\in [1,t]$ we get $x_0\in L(\o_i)$.

We prove Claim 2 by induction on $t$.
It is obvious if $t=0$. 
Since $x\in\Bb(\o_i)$ we have $\tilde e_j(x)\in\Bb(\o_i)+q^{-1}L(\o_i)$
or $\tilde e_j(x)\in q^{-1}L(\o_i)$.
In the second case we get, using Claim 1, $x_r\in q^{-1}L(\o_i)$
for all $r\in [1,t]$. Thus $x_0\in x+q^{-1}L(\o_i)$ and we are done.
Consider now the first case.
Using the induction hypothesis for $t-1$ we get an integer
$r_0\in [0, t-1]$ such that $x_{r_0+1}\in\Bb(\o_i)+q^{-1}L(\o_i)$ and
$x_{r+1}\in q^{-1}L(\o_i)$ for all $r\in[0,t-1]\setminus\{r_0\}$.
Thus $\tilde e_j(x)\in\tilde f_j^{r_0}(x_{r_0+1})+q^{-1}L(\o_i)$.
Hence
$$x\in\tilde f_j\tilde e_j(x)+q^{-1}L(\o_i)=
\tilde f_j^{r_0+1}(x_{r_0+1})+q^{-1}L(\o_i).$$
Hence, necessarily $x_0\in q^{-1}L(\o_i)$.
We are done.

We prove Claim 3.
To simplify we set $x_r=y_s=0$ if $r>t$, $s>u$.
For any $a,b\geq 0$ we have, see \cite{13, Proposition 19.1.3},
$$\la\fb_j^{(a)}(x_r)\,|\,\fb_j^{(b)}(y_s)\ra=
\delta_{a,b}\delta_{r+1,s}C_{a,s}\la x_r\,|\,y_s\ra,$$
where
$$C_{a,s}=q^{a^2-a(\mu+s\a_j,\a_j)}
\left[\matrix (\mu+s\a_j,\a_j)\cr a\endmatrix\right].$$
From Claim 1 we have $y_s\in L(\o_i)_{\mu+s\a_j}$ for all $s$.
Since $\eb_j(y_s)=0$ we have $\fb_j^{(s-1)}(y_s)=\tilde f_j^{s-1}(y_s)$.
In particular, $\fb_j^{(s-1)}(y_s)\in L(\o_i)_{\mu+\a_j}$.
Thus $\la x\,|\,\fb_j^{(s-1)}(y_s)\ra\in\ZZ[q^{-1}]$.
On the other hand,
$$\la x\,|\,\fb_j^{(s-1)}(y_s)\ra=
\sum_{r=0}^t\la\fb_j^{(r)}(x_r)\,|\,\fb_j^{(s-1)}(y_s)\ra=
C_{s-1,s}\la x_{s-1}\,|\,y_s\ra.$$
Now
$$C_{s-1,s}=q^{(s-1)^2-(s-1)(\mu+s\a_j,\a_j)}
\left[\matrix (\mu+s\a_j,\a_j)\cr s-1\endmatrix\right]
\in 1+q^{-1}\ZZ[q^{-1}].$$
Thus $\la x_{s-1}\,|\,y_s\ra\in\ZZ[q^{-1}]$ for all $s\geq 1$.
If $s\neq r+1$ then $\la x_r\,|\,y_s\ra=0$ 
since $x_r$, $y_s$ have different weights.
\quad\qed
\enddemo

\noindent
Here is the main result of the paper.

\proclaim{Theorem B}
\roster
\item The sets $\Bc_{T_\l}$, $\Bc'_{T_\l}$
are signed basis of $\Wb_{T_\l}$, $\Wb'_{T_\l}$.
Moreover, for any $\bb,\bb'\in\Bc_{T_\l}$ 
(resp. $\bb,\bb'\in\Bc'_{T_\l}$) we have
$$\partial(\bb|\bb')\in\delta_{\bb,\bb'}+q^{-1}\ZZ[q^{-1}]\quad
(\roman{resp.}\ \partial(\bb|\bb')'\in\delta_{\bb,\bb'}+q^{-1}\ZZ[[q^{-1}]]).$$

\item The signed basis $\Bc_{T_\l}$, $\Bc'_{T_\l}$
are dual with respect to the pairing $\partial(\,||\,)$
in the following sense :
for all $\bb\in\Bc_{T_\l}$ there is a unique element
$\bb'\in\Bc'_{T_\l}$ such that $\partial(\bb||\bb')=1$, and
we have $\partial(\bb||\bb'')=0$ whenever $\bb''\neq\pm\bb'$.
\endroster
\endproclaim

\demo{Proof}
Fix a decomposition $\l=\sum_{k=1}^\ell\o_{i_k}$.
We set $v_{i_1,...,i_\ell}=v_{\o_1}\otimes\cdots\otimes v_{\o_\ell}$.
By \cite{9, \S 8} the $\Ub$-submodule
$$\Nb:=\Ub\cdot(\Rb^{T_\l}\otimes v_{i_1,...,i_\ell})\subset\Otimes_k V(\o_{i_k})$$
admits a unique involution $c^\nor$ such that
$$c^\nor(v_{i_1,...,i_\ell})=v_{i_1,...,i_\ell},\quad
c^\nor(u\cdot v_{i_1,...,i_\ell})=\bar u\cdot
v_{i_1,...,i_\ell},\quad\forall u\in\Ub.$$
Set $\Rb^{T_{\o_{i_k}}}=\ZZ[z_k^{\pm 1}]$ for all $k$.
Set also 
$$\hat\Rb_{i_1/.../i_\ell}=\ZZ[[q^{-1},z_{k+1}/z_{k};k]]
\otimes_{\ZZ[q^{-1},z_{k+1}/z_{k};k]}\Rb^{\tilde T_\l},$$
where $k$ takes all possible values in $[1,\ell]$.
The tensor product 
$$\hat{\Otimes}_k\Wb_{T_{\o_{i_k}}}:=
\hat\Rb_{i_1/.../i_\ell}\otimes_{\Rb^{\tilde T_\l}}
\Otimes_k\Wb_{T_{\o_{i_k}}}$$
is endowed with the unique pairing of 
$\hat\Rb_{i_1/.../i_\ell}$-modules such that
$$(\Otimes_k x_k\,|\,\Otimes_k y_k)_{i_1/.../i_\ell}=
\Otimes_k (x_k\,|\,y_k),\quad\forall x_k, y_k\in\Wb_{T_{\o_{i_k}}},$$
and the pairing $(\,|\,)$ on each factor is as in 6.3.
As in 7.1 let $\partial(\,|\,)_{i_1/.../i_\ell}$
be the corresponding pairing
$$\hat{\Otimes}_k\Wb_{T_{\o_{i_k}}}\times
\hat{\Otimes}_k\Wb_{T_{\o_{i_k}}}\to\ZZ((q^{-1})).$$
We have an isomorphism of $\Ub$-modules
$\Wb_{T_{\o_i}}\simeq V(\o_i)$ such that $1_{\o_i}\mapsto v_{\o_i}$, 
see Theorem 7.2.A.1.
By Theorem 7.2.A.2 we have
$$\pm\Otimes_k\Bb(\o_{i_k})+\hat{\Otimes}_k q^{-1}L(\o_{i_k})=
\{\bb\in\hat{\Otimes}_k\Wb_{T_{\o_{i_k}}}\,|\,
\partial(\bb|\bb)_{i_1/.../i_\ell}=1+q^{-1}\ZZ[[q^{-1}]]\}.$$
Thus, by \cite{9, Theorem 8.5 and Proposition 8.6} the set
$$\{\bb\in\Nb\,|\,c^\nor(\bb)=\bb,\,
\partial(\bb|\bb)_{i_1/.../i_\ell}=1+q^{-1}\ZZ[q^{-1}]\}$$
is a signed basis of $\Nb$.
Iterating $(\ell-1)$-times the map in Lemma 6.2.2 
we get an isomorphism of $\Ub$-modules
$$\varpi_{i_1/.../i_\ell}\,:\,
\hat{\Otimes}_k\Wb_{T_{\o_{i_k}}}
\simto
\hat\Rb_{i_1/.../i_\ell}\otimes_{\Rb^{\tilde T_\l}}
\Wb_{T_\l}
,\quad
v_{\o_1}\otimes\cdots v_{\o_\ell}
\mapsto 
1_\l
.$$
The $\Ub$-module $\Wb_{T_\l}$ is generated by 
$\Rb^{T_\l}\otimes 1_\l$ by Lemma 5.5.2.
Thus 
$\varpi_{i_1/.../i_\ell}(\Nb)=\Wb_{T_\l}$.
By \cite{7} the pairing $(\,|\,)_{i_1/.../i_\ell}$ 
still satisfies Proposition 6.3.3
with $1_\l$ instead of $v_{i_1,...,i_\ell}.$
It is easy to see that Proposition 6.3.2 also holds
in this setting.
Thus, by Proposition 6.3.4 the pairings
$(\,|\,)_{i_1/.../i_\ell}$ and 
$(\,|\,)$ on $\Nb$
and $\Wb_{T_\l}$ coincide.
Moreover, Proposition 6.1 and Corollary 5.8.2 give 
$\varpi_{i_1/.../i_\ell}c^\nor=\b_\l\varpi_{i_1/.../i_\ell}.$
Thus, 
$$\varpi_{i_1/.../i_\ell}\bigl(
\{\bb\in\Nb\,|\,c^\nor(\bb)=\bb,\,
\partial(\bb|\bb)_{i_1/.../i_\ell}=1+q^{-1}\ZZ[q^{-1}]\}
\bigr)=$$
$$=
\{\bb\in\Wb_{T_\l}\,|\,\b_\l(\bb)=\bb,\,
\partial(\bb|\bb)\in 1+q^{-1}\ZZ[q^{-1}]\}
.$$
In particular, $\Bc_{T_\l}$ is a signed basis of $\Wb_{T_\l}$ such that
$$\partial(\bb|\bb')\in\delta_{\bb,\bb'}+q^{-1}\ZZ[q^{-1}]\leqno(7.2.11)$$
for all $\bb,\bb'\in\Bc_{T_\l}$.
By Proposition 6.3.6 there is a signed basis 
$\Bc_{T_\l}^*\subset\Wb'_{T_\l}$ dual to $\Bc_{T_\l}$
with respect to the pairing $\partial(\,||\,)$.
Let $\bb^*\in\Bc^*_{T_\l}$ be the element dual to $\bb\in\Bc_{T_\l}$.
Let $E\,:\,\Wb_{T_\l}\to\Wb'_{T_\l}$ be the unique
$\AA$-modules isomorphism such that
$E(\bb)=\bb^*$ for all $\bb\in\Bc_{T_\l}$. 
Using (7.2.11) we get
$$(\kappa_*-E)\bigl(\Oplus_{\bb\in\Bc_{T_\l}}\ZZ[q^{-1}]\bb\bigr)\subset
\Oplus_{\bb\in\Bc_{T_\l}}q^{-1}\ZZ[q^{-1}]\bb^*.$$
Thus, the map
$\kappa_*\,:\,\hat\Wb_{T_\l}\to\hat\Wb'_{T_\l}$
is invertible (see also Lemma 6.2.1) and we have
$$(\kappa_*)^{-1}=\sum_{n\geq 0}(-1)^nE^{-1}\bigl((\kappa_*-E)E^{-1}\bigr)^n.$$
In particular we get
$$\partial\bigl((\kappa_*)^{-1}(\bb^*)|{\bb'}^*\bigr)
\in\delta_{\bb,\bb'}+q^{-1}\ZZ[[q^{-1}]],$$
i.e.
$$\partial(\bb^*|{\bb'}^*)'\in\delta_{\bb,\bb'}+q^{-1}\ZZ[[q^{-1}]],$$
for all $\bb,\bb'\in\Bc_{T_\l}$.
Using Proposition 6.3.5 we get also $\b'_\l(\bb^*)=\bb^*$ 
for all $\bb\in\Bc_{T_\l}$.
Thus $\Bc^*_{T_\l}\subset\Bc'_{T_\l}$.
Then, apply Proposition 7.1.1.
We are done.
\quad\qed
\enddemo

\remark{ Remarks}
\roster
\item It is easy to see that $\Xb^{T_\l}$ is a subgroup of the multiplicative
group of $\Rb^{\tilde T_\l}$ such that $\Xb^{T_\l}\,\Bc_{T_\l}=\Bc_{T_\l}$.
Thus there is a subset $\Bb^0_{T_\l}\subset\Bc_{T_\l}$ which is
an $\Rb^{\tilde T_\l}$-basis of $\Wb_{T_\l}$.
In particular for any maximal ideal $I\subset\Rb^{\tilde T_\l}$
the set $\Bb^0_{T_\l}\otimes 1$ is a $\ZZ$-basis of 
$\Wb_{T_\l}\otimes_{\Rb^{\tilde T_\l}}(\Rb^{\tilde T_\l}/I).$
If $\dagger(I)=I$ then the involution $\b_{T_\l}$ and the metric
$\partial(\,|\,)$ descend to
$\Wb_{T_\l}\otimes_{\Rb^{\tilde T_\l}}(\Rb^{\tilde T_\l}/I).$
It is not clear if $\pm\Bb^0_{T_\l}\otimes 1$ 
admits a similar characterization as $\Bc_{T_\l}$ in 7.1.

\item Probably the sets $\Bc_\l$, $\Bc'_\l$
are signed bases of $\Wb_\l$, $\Wb'_\l$.
The conjectures in \cite{9, \S 13} and the previous theorem
suggest  that Kashiwara's canonical basis of $V(\l)$ 
coincide with $\Bc_\l$, $\Bc'_\l$, up to signs.
\endroster
\endremark
\subhead 7.3\endsubhead
We do not assume any more that $\gen$ is simply laced.
Fix $i\in I$.
The fundamental module $W(\o_i)$ is as in \cite{9}.
Let $W(\o_i)[z^{\pm 1}]$ be the affinized module,
see \cite{9, \S 4.2}.
Set $d_i=\max(1,(\a_i,\a_i)/2)$.
Then $V(\o_i)$ is isomorphic to the $\Ub$-submodule
$W(\o_i)[z^{\pm d_i}]\subset W(\o_i)[z^{\pm 1}]$,
see \cite{9, Theorem 5.15.$(viii)$}.
Set $z_i=z^{d_i}\Id\,:\,V(\o_i)\to V(\o_i)$.
Let us mention the following fact,
which is not used in the paper.

\proclaim{Proposition}
For any $\gen$ (not necessarily simply laced)
and any $i\in I$ there is a unique pairing of $\AA$-modules
$\la\,|\,\ra\,:\,V(\o_i)\times V(\o_i)\to\AA$ such that
$$\la z_i^nv_{\o_i}|z_i^mv_{\o_i}\ra=\delta_{n,m},\quad
\la u\cdot x|y\ra=\la x|\psi(u)\cdot y\ra.$$
This pairing is perfect and symmetric.
\endproclaim

\demo{Proof}
The proof is similar to the proof of \cite{13, Proposition 19.1.2}.
By \cite{9, Proposition 5.14.$(iii)$} 
the space $V(\o_i)_{\hat\mu}$ is finite-dimensional
for any $\hat\mu\in\hat P$.
Set 
$$V(\o_i)^*=\bigoplus_{\hat\mu\in\hat P}\Hom(V(\o_i)_{\hat\mu},\AA).$$
Since $\psi$ is an antiautomorphism,
there is a unique $\Ub$-module structure on $V(\o_i)^*$
such that
$$(u\cdot f)(x)=f(\psi(u)\cdot x),\quad\forall u\in\Ub,\,
\forall x\in V(\o_i).\leqno(7.3.1)$$
The $\Ub$-module $V(\o_i)^*$ is endowed with the $\hat P$-grading such that
$$V(\o_i)^*_{\mu+n\delta}=\Hom(V(\o_i)_{\mu+n\delta},\AA).$$
Recall that $V(\o_i)_{\o_i}=\AA\cdot v_{\o_i}$.
Let $f_{\o_i}\in V(\o_i)^*$ be the unique 
linear form such that $f_{\o_i}(v_{\o_i})=1$, and
$f_{\o_i}(v)=0$ for all $v\in V(\o_i)_{\hat\mu}$ with $\hat\mu\neq\o_i$.
Hence, $f_{\o_i}\in V(\o_i)^*_{\o_i}$.
We must prove that there is a unique morphism of $\Ub$-modules 
$V(\o_i)\to V(\o_i)^*$ which takes $v_{\o_i}$ to $f_{\o_i}$,
and that it is invertible.
The spaces $V(\o_i)_{\hat\mu}$, $V(\o_i)^*_{\hat\mu}$ 
have the same dimension for all $\hat\mu$.
Thus the set of the weights $\mu\in P$ 
such that $V(\o_i)^*_{\mu+n\delta}\neq\{0\}$ for some $n\in\ZZ$
is contained in $\o_i-\sum_{j\in I}\NN\,\a_j$, 
since this is true for $V(\o_i)$.
Hence $f_{\o_i}$ is an extremal vector of weight $\o_i$, 
see \cite{9, Theorem 5.3}.
By the universal property of $V(\o_i)$,
there is a unique morphism of $\Ub$-modules 
$\phi\,:\,V(\o_i)\to\Ub\cdot f_{\o_i}\subseteq V(\o_i)^*$ 
which takes $v_{\o_i}$ to $f_{\o_i}$.
Moreover, we have $\phi(V(\o_i)_{\hat\mu})\subseteq V(\o_i)^*_{\hat\mu}$ 
for all $\hat\mu\in\hat P$.
Since $V(\o_i)_{\hat\mu}$ is finite dimensional it is sufficient to prove
that the map $\phi$ is injective.
Let the operator $z_i$ acts on $V(\o_i)^*$ by
$$(z_i\cdot f)(x)=f(z_i^{-1}\cdot x),
\quad\forall x\in V(\o_i),\forall f\in V(\o_i)^*.$$
Then $\phi$ commutes to $z_i$.
Since $W(\o_i)\simeq V(\o_i)/(z_i-1)V(\o_i)$,
the map $\phi$ induces a non-zero morphism of 
$\Ub$-modules $W(\o_i)\to W(\o_i)$.
It is injective since $W(\o_i)$ is simple.
Thus $\phi$ is injective.

The pairing is symmetric because it is
unique and $\psi^2=\Id$.
\quad\qed
\enddemo
\head 8.  Example\endhead
We assume that $\Pi=(\{1\},\emptyset)$. 
We set $\l=\ell\o_1$, $\a=a\a_1$.
To simplify we omit the subscripts 1 and we set
$Q_{\ell a}=Q_{\l\a}$, $d_{\ell a}=d_{\l\a}$, etc.
Set
$$\tilde x_{r}^+=\sum_{a'=a+1}(-1)^{\ell-a'}
q^{-2r-a'}(\Wedge_\Vc^{-r-\ell+a'}\boxtimes\Wedge_\Vc^{r+\ell-a})
\otimes\Wedge_{\Wc}^{-1}\otimes 1_{\ell aa'},$$
$$\tilde x_{r}^-=\sum_{a'=a-1}(-1)^{a'}
q^{-2r+a'}(\Wedge_\Vc^{r-a'}\boxtimes\Wedge_\Vc^{-r+a})
\otimes 1_{\ell aa'}.$$
We have $Q_{\ell a}\neq\emptyset$ if and only if $0\leq a\leq\ell$. 
More precisely $F_{\ell a}$ is smooth and isomorphic the Grassmanian
of $a$-dimensional subspaces in $\CC^\ell$,
and $Q_{\ell a}=T^*F_{\ell a}$.
An element in $Q_{\ell a}$ may be viewed as a couple $(V,u)$,
where $V\subseteq\CC^\ell$ is a $a$-dimensional subspace,
and $u\in\End(\CC^\ell)$ is a nilpotent map such that
$\Im u\subseteq V\subseteq \Ker u$.
The element $z\in\CC^*$ acts on $T^*F_{\ell a}$
by multiplication by the scalar $z^2$ along the fibers.
The automorphism $\o\,:\,Q_{\ell a}\to Q_{\ell, \ell-a}$
takes the couple $(V,u)$ to $(V^\perp,{}^t u)$, where
$V^\perp\subseteq\CC^\ell$ is the subspace orthogonal to $V$,
with respect to the canonical scalar product on $\CC^\ell$.
Let $\Ec'_{a}$ be the tautological rank $a$ vector bundle on $T^*F_{\ell a}$.
Let $\Ec_{a}$ be the restriction of $\Ec'_{a}$ to $F_{\ell a}$.
Set $\Qc'_{a}=\Wc_{a}/\Ec'_{a}$,
$\Qc_{a}=\Wc_{a}/\Ec_{a}$.
Set $\Ec=\Oplus_a\Ec_{a}$, $\Ec'=\Oplus_a\Ec'_{a}$, etc.
We have 
$$\Vc=q\Ec',\quad
\tilde x_r^+=\Wedge^{-\ell}_{\Ec'}\
x_r^+\Wedge_\Wc^{-\ell}\Wedge^{\ell}_{\Ec'},\quad
\tilde x_r^-=\Wedge^{-\ell}_{\Ec'}
x_r^-\Wedge_\Wc^\ell\Wedge^{\ell}_{\Ec'}.$$
Hereafter we omit the operators $\kappa_*$ and $\otimes$. 
The following lemma is immediate, see for instance \cite{22}.

\proclaim{Lemma}
\roster
\item We have $\tilde x^+_0(1_{\ell a})=[\ell-a+1]1_{\ell,a-1},$
$\tilde x^-_0(1_{\ell a})=[a+1]1_{\ell,a+1}.$

\item We have 
$1_{\ell a}=\sum_{i=0}^{a(\ell-a)}(-1)^i q^{-2i} 
\wedge^i\bigl({\Qc'}_{a}{\Ec'_a}^*\bigr).$
\endroster
\endproclaim

\noindent
A direct computation gives

\proclaim{Proposition} 
\roster
\item We have
$\Bc'_1=\Bc_1=\pm\{x1_{10},x\Wedge_\Wc^{-1}1_{11}; x\in\Xb_1\}.$

\item We have
$\Bc'_2=\pm\{x1'_{20},x\Wedge_\Wc^{-2}1'_{21},
q^{-1}x\Wedge_\Wc^{-2}\Qc'_{21},x\Wedge_\Wc^{-4}1'_{22}; 
x\in\Wedge_{\Ec'_2}^2\otimes\Xb_2\},$

\noindent
$\Bc_2=\pm\{x1_{20},x\Wedge_\Wc^{-2}1_{21},q^{-1}x\Wedge_\Wc^{-2}\Ec_{21},
x\Wedge_\Wc^{-4}1_{22}; x\in\Wedge_{\Ec_2}^2\otimes\Xb_2\}.$
\endroster
\endproclaim

\head Appendix \endhead
Let us check that $x(w_0*\a)=x(\a)$ and $y(w_0*\a)=y(\a),$  for all $\a\in Q,$
where $x\,:\, Q\to\ZZ$ and $y\,:\, Q\to \ZZ/2\ZZ$ are the quadratic maps
which satisfy (5.7.1).  We give the proof for the map $x$, the case of the map
$y$ is left to the reader. Set
$$\a=\sum_ja_j\a_j,\quad \l=\sum_j\ell_j\o_j,\quad\nu=\l-w_0(\l)=
\sum_jk_j\a_j.$$
Write $x(\a)=Q(\a)+L(\a)+a,$ where $Q$ is 
a quadratic form, $L$ is a linear form and $a\in\ZZ.$ Put
$$Q(\a)=\sum_{i,j}q_{ij}a_ia_j,\quad
L(\a)=\sum_jb_ja_j.$$
A direct computation gives $g_{i;\a}=c_i+\sum_jc_{ij}a_j$ with 
$$c_i=-1+(2-\cb)\ell_i+(2\cb-1)k_i,
\quad c_{ij}=(\cb-2)a_{ij}+\delta_{ij}(2-2\cb).$$
Using the  relation $x(\a+\a_i)-x(\a)=\cb k_i-g_{i;\a},$
we get
$$ b_i=\cb k_i-c_i+{1\over 2}c_{ii}=(\cb -2)\ell_i+(1-\cb)k_i,$$
$$q_{ij}=-{1\over 2}c_{ij}=(1-{\cb\over 2})a_{ij}+\delta_{ij}(\cb-1),
\leqno(1)$$
that is 
$$Q(\a)=(\cb-1)|\a|^2+(1-{\cb\over 2})(\a,\a), \quad L(\a)=(\cb-2)(\l,\a)
+(1-\cb)\sum_jk_ja_j.\leqno(2)$$
The identity $Q(\nu)+L(\nu)+a=x(\nu)=(\cb -1)|\nu|^2-\cb (\l,\l),$ gives
$$a=-\cb(\l,\l)+({\cb\over 2}-1)(\nu,\nu)-(\cb-2)(\l,\nu)
+(\cb-1)|\nu|^2= (\cb-1)|\nu|^2-\cb(\l,\l),$$
since  $2(\l,\nu)=(\nu,\nu).$ In particular 
$$a=x(0)=x(\nu).$$
We want now to check that $x(w_0*\a)=x(\a).$ 
Since  $x(\nu)=x(0),$ we have $Q(\nu)+L(\nu)=0.$ Moreover (2) gives
$Q(\a)=Q(w_0(\a)).$ Thus it is enough to prove that 
$$L(w_0(\a))-L(\a)=2\sum_{i,j}k_iq_{ij}a_{\uj}.$$
By (2), the left hand side is equal to 
$$(2-\cb)(\nu,\a)+(\cb-1)\sum_jk_j(a_j+a_{\uj}).$$
By (1), the right hand side  is equal to 
$$(2-\cb)(\nu,\a)+2(\cb-1)\sum_jk_ja_{\uj}.$$
Since $w_0(\nu)=-\nu$ we have $k_j=k_\uj.$ We are done. 

\vskip1cm
\refstyle{C}

\vskip3cm
\enddocument